\documentclass[11pt,a4paper]{article}
\usepackage[utf8]{inputenc}
\usepackage{amsmath,bm}
\usepackage{enumerate}
\usepackage{amsmath,nccmath}
\usepackage{amsthm}
\usepackage{amsfonts}
\usepackage[square, comma, sort&compress,numbers]{natbib}
\usepackage{amssymb}
\usepackage{setspace}
\usepackage{graphicx}
\usepackage{epstopdf}
\usepackage[shortlabels]{enumitem} 
\usepackage{verbatim}
\usepackage{color}
\usepackage[colorlinks=true,linktocpage=true]{hyperref}
\usepackage[marginal]{footmisc}
\usepackage[section]{placeins}
\usepackage{indentfirst}
\usepackage[left=2cm,right=2cm,top=2cm,bottom=2cm]{geometry}
\usepackage{float}
\usepackage{hyperref} 
\usepackage[doipre={DOI:~}]{uri} 
\usepackage{enumitem}
\usepackage{tocloft}
\usepackage{titlesec} 
\titleformat{\subsubsection}[runin]
  {\normalfont\normalsize\bfseries}{\thesubsubsection}{1em}{}

\author{}
\theoremstyle{plain}
\newtheorem{prop}{Proposition}[section]
\newtheorem*{prop*}{Proposition}
\newtheorem{cor}{Corollary}
\newtheorem{lem}[prop]{Lemma}
\newtheorem{thm}{Theorem}
\newtheorem*{lem*}{Lemma}
\newtheorem*{thm*}{Theorem}
\newtheorem{claim}[prop]{Claim}
\newtheorem*{claim*}{Claim}

\newtheorem{thma}{Theorem}

\newtheorem{thms}[prop]{Theorem}

\theoremstyle{definition}
\newtheorem{defi}[prop]{Definition}
\newtheorem*{defi*}{Definition}

\newtheorem*{ques*}{Question}
\newtheorem*{conj*}{Conjecture}

\theoremstyle{remark}
\newtheorem{rem}[prop]{Remark}

\allowdisplaybreaks
\usepackage{mathtools}

\DeclarePairedDelimiter\floor{\lfloor}{\rfloor}

\usepackage{authblk}

\renewcommand{\geq}{\geqslant}
\renewcommand{\leq}{\leqslant}

\definecolor{darkred}{rgb}{.56,0,0}

\makeatletter
\newcommand*{\rom}[1]{\expandafter\@slowromancap\romannumeral #1@}
\newcommand\primeitem{%
 \item[(\roman{enumi}\textquotesingle)]\def\@currentlabel{(\roman{enumi}\textquotesingle)}}
\newcommand\pprimeitem{%
 \item[\arabic{enumi}\textquotesingle]\def\@currentlabel{\arabic{enumi}\textquotesingle}}
 
\newcommand{\setword}[2]{%
  \phantomsection
  #1\def\@currentlabel{\unexpanded{#1}}\label{#2}%
}

\newcommand{\address}[1]{\gdef\@address{#1}}
\newcommand{\email}[1]{\gdef\@email{\url{#1}}}
\newcommand{\@endstuff}{\par\vspace{\baselineskip}\noindent\small
\begin{tabular}{@{}l}\scshape\@address\\\textit{E-mail address:} \@email\end{tabular}}
\AtEndDocument{\@endstuff}

\makeatother

\title{{\bf Blender-producing mechanisms and a dichotomy for local dynamics near heterodimensional cycles}\\
}
\author{Dongchen Li}
\address{Shanghai Center for Mathematical Sciences, Fudan University, China}
\email{dongchenli@fudan.com.cn}

\begin{document}
\def\D{\mathrm{D}}
\def\d{\mathrm{d}}
\def\tr{\mathrm{tr}\,}
\def\diff{Di\! f\! f}
\def\eps{\varepsilon}

\def\q2{8}
\def\simple{simple hyperbolic dynamics}

\def\ind{\operatorname{ind}}
\def\arg{\operatorname{Arg}}
\def\rank{\operatorname{rank}}
\def\card{\operatorname{Card}}

\bibliographystyle{plainnat}

\maketitle

\begin{flushright}
\textit{Dedicated to Professor Dmitry Turaev on his 60th birthday}
\end{flushright}
\vspace{.3cm}
\par{}
\noindent{\bf Abstract.} 
Blenders are special hyperbolic sets  used to produce various robust dynamical phenomena which appear fragile at  first glance.
We prove for $C^r$  diffeomorphisms ($r=2,\dots,\infty,\omega$) that
blenders naturally exist (without perturbation) near non-degenerate heterodimensional cycles of coindex-1, and the existence 
is determined by arithmetic properties of moduli of topological conjugacy for diffeomorphisms with heterodimensional cycles. In particular, we obtain a dichotomy for  dynamics in any  small neighborhood $U$ of a non-degenerate heterodimensional cycle: 
either there exist   infinitely many blenders  accumulating on the cycle, forming robust heterodimensional dynamics in most cases, or
there are no orbits other than those constituting the cycle  lying entirely in $U$.




\tableofcontents

\section{Introduction}\label{sec:intro}
A blender is a hyperbolic set such that its projection to a certain ``central'' subspace has a non-empty interior, and this property persists for 
any $C^1$-close system. An immediate consequence is that non-transverse intersections with the stable or unstable set can be unremovable by $C^1$-small perturbations. This opens the way for constructing various $C^1$-robust dynamical phenomena, including those which appear fragile at  first glance.
As the first such example, a large class of diffeomorphisms exhibiting the robust nonhyperbolic transitivity was constructed by Bonatti and D\'iaz in their discovery paper \citep{BD96} of blenders. 
Subsequently, various non-trivial results involving blenders have been obtained: the persistence of heterodimensional cycles \citep{BD08,BDK12,LT21},  abundance of $C^1$-robust homoclinic tangencies \citep{BDP22,BD12,LLST22}, typicality of the Newhouse phenomenon in the space of families of diffeomorphisms \citep{Be16}, robustly fast growth of the number of periodic orbits \citep{AST17,AST21,Be21},  robust existence of nonhyperbolic ergodic measures \citep{BBD16}, $C^1$ density of stable ergodicity \citep{RRTU11,ACW21}, robust transitivity in Hamiltonian dynamics \citep{NP12},  persistence of homoclinics to saddle-center periodic points \citep{LT23}, and robust bifurcations in complex dynamics \citep{Bi19,Du17,Ta17}, among others. 

Essentially, there are  two strategies for obtaining a blender. One can  either  explicitly construct blenders   from (relatively) specific  systems, see e.g. \citep{Be16,NP12,Be21}, or identify conditions under which blenders emerge near particular  dynamical objects. In the latter case, blenders are mostly found near heterodimensional cycles, except that, in \citep{MS12,ACW21},  blenders  arise from  horseshoes with sufficiently large  fractal dimension, while in \citep{LT23} blenders are built near one-dimensional whiskered tori with homoclinics. In this paper, we   investigate the creation of  blenders near heterodimensional cycles. The first major contribution in this direction was made by D\'iaz in \citep{Di95a} where he implicitly constructed a blender. Later, Bonatti and D\'iaz developed \citep{BD08} a general theory in the $C^1$ topology, establishing that  blenders can appear after a $C^1$-small perturbation which unfolds a heterodimensional cycle of coindex 1. A conservative version  of this result was also obtained in \citep{RRTU10}. More recently, together with Turaev, the author showed in \citep{LT21} that, for every $C^r$  diffeomorphisms $(r=2,\dots,\infty,\omega)$ having a coindex-1 heterodimensional cycle, there exist $C^r$-small perturbations which give rise to blenders. The method employed there is different from \citep{BD08} and, in particular, we identified certain classes of heterodimensional cycles which give rise to blenders without being destroyed.

The present paper serves as a continuation of the blender-related investigation initiated in \citep{LT21}, completing the analysis of the case where  at least one of the saddles in a heterodimensional cycle has a nonreal central multiplier.  The main goal is to provide a comprehensive description (Theorem~\ref{thm:main}) of the relation between the natural (i.e., non-perturbative) existence of blenders and  arithmetic properties of  moduli of topological conjugacy of the $C^r$ diffeomorphisms\footnote{The results in this paper can be translated without effort to the setting of flows on manifolds of dimension four or higher, see Remark \ref{rem:flow}.} having heterodimensional cycles of coindex 1.  The result  characterizes the blender-producing mechanisms for coindex-1 heterodimensional dynamics in full generality. As a consequential result, we establish a dichotomy (Corollary~\ref{cor:limit} and Theorem~\ref{thm:rational})  for the local dynamics near  non-degenerate heterodimensional cycles, depending on the existence/non-existence of rational dependencies between the moduli: one has either  rich dynamics, that is, the cycle is a limit of infinitely many blenders, which in most cases form robust heterodimensional dynamics; or simple hyperbolic dynamics, that is,   there are only finitely many locally maximal and compact hyperbolic sets in any sufficiently small neighborhood of the cycle and no heterodimensional connections (except for the cycle itself). Moreover, such dichotomy can extend to nearby systems (Corollary~\ref{cor:dich}). We also obtain several perturbative results (Theorem~\ref{thm:unfold}) which can be used as a standard tool for creating blender-involved robust heterodimensional dynamics in the $C^r$ topology, giving possibilities to improve the regularity of many $C^1$ results  based on the Bonatti-D\'iaz construction in \citep{BD08}.

 
\subsection{Blenders and moduli of topological conjugacy}\label{sec:moduli}
There are many definitions of blenders, see e.g. \citep{BD96,BD12,BBD16,BDV,BKR14,NP12}. The blenders found in this paper are the same as in \citep{LT21}, which we call {\em standard}. They carry more dynamical meaning than the original geometrically defined blenders \citep{BD96,BDV,BD08}. Roughly speaking, a standard blender is a zero-dimensional partially hyperbolic set obtained from a Markov partition with finite (but arbitrarily large) number of elements, such that the  iterated function system which arises from the partition satisfies the so-called covering property (see Proposition \ref{prop:blender}). It is a generalization of the {\em blender-horseshoe} defined in \citep{BD12}, whose corresponding Markov partition had exactly two elements. We refer interested readers to \citep{LT21,Li24}, where a detailed construction  of standard blenders is given.  Here for stating our results, it suffices to use a version of the geometric definition \cite[Definition 6.11]{BDV}. 

In the rest of the paper, let $\mathcal{M}$ be a manifold of dimension at least three and let $\diff ^r(\mathcal{M})$ be the space of $C^r$ diffeomorphisms of $\mathcal{M}$.

\begin{defi}[Blenders]\label{defi:blender}
A basic (i.e., compact, transitive, locally maximal) hyperbolic set $\Lambda$ of a diffeomorphism $f\in\diff ^1(\mathcal{M})$ is called a {\em center-stable (cs) blender} if there exists a $C^1$-open set $\mathcal{D}$ of  $d^{ss}$-dimensional discs (embedded copies of $\mathbb{R}^{d^{ss}}$) with $d^{ss}=\dim(W^s(\Lambda))-1$, such that for every system $g$ which is $C^1$-close to $f$, for the hyperbolic continuation $\Lambda_g$ of the basic set $\Lambda$, the set $W^u(\Lambda_g)$ intersects every element from $\mathcal{D}$; a {\em center-unstable (cu) blender} is a cs-blender of $f^{-1}$.
\end{defi}

 Recall that the index of a transitive hyperbolic set is the rank of its unstable bundle. We denote it by $\ind(\cdot)$. Following \citep{LT21}, we use
\begin{defi}[Heterodimensional dynamics]\label{def:hdc}
We say that a diffeomorphism $f\in\diff ^1(\mathcal{M})$ has {\em heterodimensional dynamics} if there are two basic  hyperbolic sets $\Lambda_1$ and $\Lambda_2$ such that they have different indices and their invariant manifolds intersect cyclically, i.e., $W^u(\Lambda_1)\cap W^s(\Lambda_2)\neq \emptyset$ and $W^u(\Lambda_2)\cap W^s(\Lambda_1)\neq \emptyset$. The quantity $|\ind(\Lambda_1)- \ind(\Lambda_2)|$ is called the coindex of the heterodimensional dynamics.
\end{defi}
 In literature, the above  heterodimensional dynamics is called a `heterodimensional cycle'. However, we reserve the word `cycle' for the simplest version of such dynamics, which is  the precise object studied in the paper.
\begin{defi}[Heterodimensional cycles]\label{def:hdc}
A {\em heterodimensional cycle} of $f$ is a closed invariant set consisting of four orbits: two hyperbolic periodic orbits $L_1$ and $L_2$ with $\ind(L_1)\neq \ind(L_2)$, and two heteroclinic orbits, one from $W^u(L_1)\cap W^s(L_2)$ and the other from $W^u(L_2)\cap W^s(L_1)$. 
\end{defi}

\begin{figure}[!h]
\begin{center}
\includegraphics[width=0.55\textwidth]{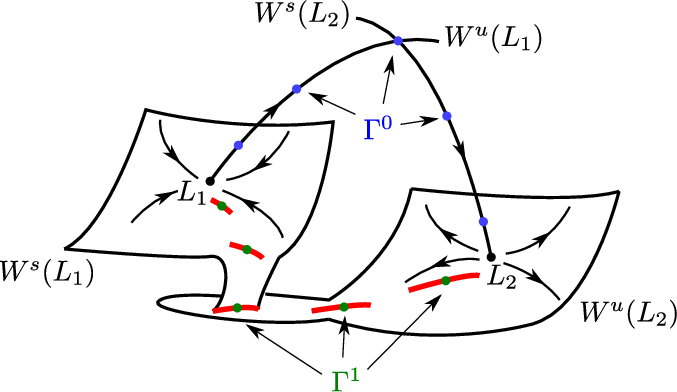}
\end{center}
\caption{\citep[Figure 1]{LT21}. A heterodimensional cycle for a three-dimensional diffeomorphism. It consists of two fixed points $L_1$ and $L_2$, a fragile heteroclinic orbit $\Gamma^0$ (blue dots) in the non-transverse intersection $W^u(L_1)\cap W^s(L_2)$, and a robust heteroclinic orbit $\Gamma^1$ (green dots) in the transverse intersection $W^s(L_1)\cap W^u(L_2)$ (red curves).}
\label{fig:hdc0}
\end{figure}

See Figure \ref{fig:hdc0} for an illustration. We will only study coindex-1 heterodimensional cycles, i.e., those with $|\ind(L_1)-\ind(L_2)|=1$, so for brevity we omit the word `coindex-1' in the remaining paper. Note that we lose no generality by considering only cycles associated with periodic orbits. Indeed, since hyperbolic periodic orbits are dense in a hyperbolic basic set and their stable/unstable invariant manifolds approximate the stable/unstable leaves of the hyperbolic set (see e.g. \cite[Theorem 6.4.15]{KH95}), one can always obtain heterodimensional cycles from heterodimensional dynamics by an arbitrarily $C^r$-small perturbation.

We further consider heterodimensional cycles that fulfill certain non-degeneracy conditions introduced in \citep{LT21} (which are also presented in Appendix~\ref{sec:nd} for completeness). All these conditions can be achieved by an arbitrarily $C^r$-small perturbation and are $C^r$-open in the set of diffeomorphisms having a heterodimensional cycle.

The dynamics near a non-degenerate heterodimensional cycle is largely determined by the central multipliers\footnote{The multipliers of a periodic orbit $L$ are the eigenvalues of the derivative of the period map $f^{per(L)}$ evaluated at any point of $L$.} of the two involved periodic orbits $L_1$ and $L_2$.  Let us assume $\ind(L_1)+1=\ind(L_2)$, and denote $d_1:=\ind(L_1)$ and $d_2:=\ind(L_2)$. When the multipliers of $L_1$ are in general position, there is either one real multiplier, say $\lambda$, or a pair of conjugate complex multipliers, say $\lambda$ and $\lambda^*$, closest to the unit circle from inside, which we call {\em center-stable}. Similarly, for $L_2$, there is either one real multiplier, say $\gamma$, or a pair of conjugate complex multipliers, say $\gamma$ and $\gamma^*$, closest to the unit circle from outside, which we call {\em center-unstable}. See Section \ref{sec:subcases} for a detailed discussion.  It is well-known \citep{NPT83,St82} that 
the ratio
\begin{equation*}\label{eq:theta}
\theta:=-\dfrac{\ln|\lambda|}{\ln|\gamma|}
\end{equation*}
and the arguments
\begin{equation}\label{eq:argu}
\omega_1 := \arg (\lambda)\in [0,\pi] \quad\mbox{and}\quad \omega_2 := \arg (\gamma) \in [0,\pi],
\end{equation}
are {\em moduli}, i.e., invariants of topological conjugacy, for diffeomorphisms having non-degenerate heterodimensional cycles. Note that  the arguments are also defined for real $\lambda$ and $\gamma$, in which case one has $\omega_1\in\{0,\pi\}$ and, respectively, $\omega_2\in\{0,\pi\}$. We prove
\begin{thma}\label{thm:main}
Let $f\in\diff ^r(\mathcal{M})$, with $r=2,\dots,\infty,\omega$,  have a non-degenerate heterodimensional cycle, and let $U$ be any small neighborhood of the cycle.
\begin{itemize}[nosep]
\item If $\omega_1/2\pi$ and $\omega_2/2\pi$ are both irrational, then there exist both a cs-blender and a cu-blender in $U$. The blenders have different indices and their stable and unstable sets intersect cyclically and $C^1$-robustly;
\item If exactly one of $\omega_1/2\pi$ and $\omega_2/2\pi$ is irrational, then there exists at least one blender
\footnote{This blender, and also the one in the third item, can be center-stable or center-unstable. This depends on  $f$, and will be make precise in Propositions \ref{prop:blendersf} and \ref{prop:blenderdf}}
 in $U$;
\item If  $\theta$ is irrational with at least one of  $\lambda$ and  $ \gamma$ being nonreal, and $\omega_1/2\pi$ and  $\omega_2/2\pi$ do not belong to the set
\begin{equation*}\label{eq:set}
\left\{\dfrac{p}{q}\mid p,q\in\mathbb{N} \;\mbox{are coprime with } 3\leq q \leq 8\right\},
\end{equation*}
then there exists at least one blender in $U$.
\end{itemize}
In each of the above cases, the blender is center-stable if and only if  it has index $d_1$, and center-unstable if and only if it has index $d_2$.
\end{thma}
\begin{rem}
This theorem significantly strengthens the previous results \citep[Propositions 6.2 and 6.7]{LT21} in terms of the required conditions, where only the following three subcases are considered: (1) $\theta,\omega_1/2\pi,1$ are rationally independent, (2) $\theta,\omega_2/2\pi,1$ are rationally independent, and (3) $\theta,\omega_1/2\pi,\omega/2\pi,1$ are rationally independent.
\end{rem}

By definition, if $\omega_1/2\pi$ is irrational, the multiplier $\lambda$ must be nonreal, and, if it is rational, $\lambda$ can be either real or nonreal; similarly for $\omega_2/2\pi$ and $\gamma$. Hence, the result concerns all the possibilities of the modulus values when at least one of $\lambda$ and $\gamma$ is nonreal, except for the possibility that all $\omega_{1}/2\pi,\omega_{2}/2\pi$ and $\theta$ are rational, considered in Theorem \ref{thm:rational} below. Theorem \ref{thm:main} is our main result and it reflects the essential relation between values of the moduli and the existence of blenders.  However, it is too general for practical use in the sense that not much information of the blenders is provided. In fact, this theorem  follows from several much finer results collected in Theorems \ref{thm:sf} and \ref{thm:df} in Section \ref{sec:subcases}, where the homoclinic relations between the found blenders and  periodic orbits involved in the cycle are established. We stress that Theorem \ref{thm:main} has a non-perturbative nature. This is the main difference from the Bonatti-D\'iaz construction in \citep{BD08}, where the cycle is destroyed to create blenders.

The case where both $\lambda$ and $\gamma$ are real   is fully studied in \citep{LT21}, where such cycles are further classified to types I and II, according to whether an orientation-like property is satisfied by the central dynamics along the non-transverse heteroclinic orbit, see \citep[Section 2.5]{LT21}. 
\begin{thm*}[{\citep[Theorem 1]{LT21}}]\label{thm:type1}
For any non-degenerate type-I cycle and any neighborhood $U$ of it, if $\theta$ is irrational, then there exists a blender in $U$. 
\end{thm*}

Note that a cycle being type-II in particular means that both central multipliers $\lambda$ and $\gamma$ are  real. In summary, we obtain
\begin{cor}\label{cor:limit}
For any non-degenerate heterodimensional cycle which is not type-II, it is accumulated by infinitely many blenders if
\begin{itemize}[nosep]
\item either $\dfrac{\omega_i}{2\pi}$ is irrational for some $i\in \{1,2\}$,
\item or $\theta$ is irrational, and $\dfrac{\omega_i}{2\pi}=\dfrac{p_i}{q_i}$ for  $i=1,2$, where $p_i,q_i\in\mathbb{Z}$ are coprime  with $q_i> 8$ .
\end{itemize}
\end{cor}

We remark that the infinite blenders in Corollary \ref{cor:limit} form a ``nearly-affine blender system'', as defined in \citep{Li24}.  The union of a properly chosen  finite collection of  those blenders  carries  dynamical properties which normally do not hold for a single blender, which are discussed in details in \citep{Li24}.
 In a forthcoming paper \citep{LT24}, we also show that these blenders can  be  embedded into a family to form a parablender in the sense of \citep{Be16}, giving rise to robust coexistence of infinitely many  sinks for generic families of diffeomorphisms.

\subsection{Simple hyperbolic dynamics}\label{sec:simple}
In the last remaining case where  $\theta$, $\omega_1/2\pi$ and $\omega_2/2\pi$ are all rational (including the possibilities of  $\lambda$ and/or $\gamma$ being real), there are generally no blenders near the heterodimensional cycle. To state the results precisely, let us  introduce some necessary notions.

Denote the heterodimensional cycle of $f$ by $\Gamma$ , and the two involved heteroclinic orbits  by $\Gamma^0\subset W^u(L_1)\cap W^s(L_2)$ and $ \Gamma^1\subset W^u(L_2)\cap W^s(L_1)$, where we assumed that $\ind(L_1)+1=\ind(L_2)$. Recall that the hyperbolic periodic orbits admit unique continuations for all diffeomorphisms $C^1$-close to $f$. Thus, for any small neighborhood $U$ of $\Gamma$, we can find  a codimension-1 subspace $\mathcal{H}_f \subset \diff ^r(\mathcal{M})$ such that each $g\in \mathcal{H}_f$ has a cycle $\Gamma_g\subset U$, involving the continuations of $L_1$ and $L_2$, and heteroclinic orbits close to $\Gamma^{0}$ and $\Gamma^1$. Since $\theta$ depends continuously on diffeomorphisms, inside $\mathcal{H}_f$ there is a codimension-2 subspace $\mathcal{H}_{f,\theta}=\mathcal{H}_f\cap \{\theta=const\}$, where all diffeomorphisms have the same $\theta$ values as $f$. Similarly, one can define two codimension-3 subspaces $\mathcal{H}_{f,\theta,\omega_1}$ and $\mathcal{H}_{f,\theta,\omega_2}$ and one codimension-4 subspace $\mathcal{H}_{f,\theta,\omega_1,\omega_2}$ by intersecting $\mathcal{H}_{f,\theta}$ with $\{\omega_1=const\}$, $\{\omega_2=const\}$ and both of them, respectively.

We define the splitting parameter $\mu$ for the heterodimensional cycle as the continuous functional of a  neighborhood $\mathcal{U}\subset \diff ^r(\mathcal{M})$ of $f$ which  measures the distance between $W^u(L_1)$ and $W^s(L_2)$ in a certain neighborhood (see Remark \ref{rem:mu0}). Then, $\mu(g)=0$ for $g\in\mathcal{U}$ corresponds to  the codimension-1 subspace $\mathcal{H}_f$. For a one-parameter family $\{f_\eps\}$ of $C^r$ diffeomorphisms with $f_0=f$, we call it a {\em generic unfolding family} (when there is one parameter) or a {\em rank-1 unfolding family} (when there is more than one parameter) if
\begin{equation}\label{eq:rank}
\rank\left.\dfrac{d\mu(\eps)}{d\eps}\right|_{\eps=0}= 1.
\end{equation}
Let $\mathcal{N}$ be the set of orbits which lie entirely in $U$. 

\begin{defi}[Simple hyperbolic dynamics]\label{defi:simple}
 A generic one-parameter unfolding family $\{f_\eps\}$ is said to have {\em\simple} if
\begin{itemize}[nosep]
\item at $\eps=0$, $\mathcal{N}$ is the union of $L_1$, $L_2$, $\Gamma^0$, and the orbits of transverse intersection of $W^u(L_2)$ with $W^s(L_1)$ near $\Gamma^1$;
\item at all sufficiently small $\eps\neq 0$, $\mathcal{N}$ is comprised 
\begin{itemize}[nosep]
\item either by $L_2$, an index-$d_1$ uniformly-hyperbolic compact set $\Lambda_1$ containing $L_1$, and transverse heteroclinic
orbits in $W^u(L_2)\cap W^s(\Lambda_1)$, 
\item  or  by $L_1$, an index-$d_2$ uniformly-hyperbolic compact set $\Lambda_2$ containing $L_2$, and transverse heteroclinic
orbits in $W^u(\Lambda_2)\cap W^s(L_1)$.
\end{itemize}
\end{itemize}
In particular, there exist no blenders at $\eps=0$ and   no  heterodimensional dynamics at $\eps\neq 0$.
\end{defi}

\begin{thma}\label{thm:rational}
Let $f\in \diff ^r(\mathcal{M})$, with $r=2,\dots,\infty,\omega$,  have a  heterodimensional cycle with all $\theta$, $\omega_1/2\pi$ and $\omega_2/2\pi$ being rational. There exists a set $\mathcal{H}'_f\subset \mathcal{H}_{f}$ such that any generic one-parameter unfolding family $\{g_\eps\}$ with $g_0\in \mathcal{H}'_{f}$ has only {\simple} defined above.
Moreover, the set $\mathcal{H}'_f$ satisfies that
\begin{itemize}[nosep]
\item if  $\lambda$ and $\gamma$ are both real, then $\mathcal{H}'_f$ is  open and dense in $ \mathcal{H}_{f,\theta}$,
\item if $\lambda$ is nonreal and $\gamma$ is real, then $\mathcal{H}'_f$ is  open and dense in $ \mathcal{H}_{f,\theta,\omega_1}$,

\item if $\lambda$ is real and $\gamma$ is nonreal, then $\mathcal{H}'_f$ is  open and dense in $\mathcal{H}_{f,\theta,\omega_2}$, and

\item if  $\lambda$ and $\gamma$ are both nonreal,  then $\mathcal{H}'_f$ is  open and dense in $\mathcal{H}_{f,\theta,\omega_1,\omega_2}$. 
\end{itemize}
\end{thma}

The first item in Theorem~\ref{thm:rational}  is  established in \citep[Theorem 6]{LT21}, and the proof of the remaining part 
 is given in Section \ref{sec:rational}. 
The result implies that the blenders which exist for the moduli values in Corollary  \ref{cor:limit} depart from the heterodimensional cycle immediately once we change the moduli to rational values within $\mathcal{H}_f$.


\subsection{Robust heterodimensional dynamics arisen from perturbations}
A direct consequence of the existence of blenders near a heterodimensional cycle is the emergence of robust heterodimensional dynamics, defined as
\begin{defi}[Robust heterodimensional dynamics]\label{defi:hdd_robust}
We say that a diffeomorphism $f$ exhibits {\em $C^1$-robust heterodimensional dynamics} if it has heterodimensional dynamics involving two  hyperbolic basic sets $\Lambda_1$ and $\Lambda_2$, where at least one of them is non-trivial, and there exists a $C^1$-neighborhood $\mathcal{U}$ of $f$ such that every $g\in \mathcal{U}$ has heterodimensional dynamics involving the continuations of $\Lambda_1$ and $\Lambda_2$.
\end{defi}
It is proven that robust heterodimensional dynamics can always be obtained by an arbitrarily small perturbation of a heterodimensional cycle, see \citep{BD08,BDK12} for the $C^1$ topology and \citep{LT21} for the $C^r$ topology. Below we present a strengthened version of the latter result, in terms of the class of families allowed. Specifically, we consider the case where at least one of $\lambda$ and $\gamma$  is nonreal. 

Let $f$ have a heterodimensional cycle involving $L_1$ and $L_2$ with $\ind(L_1)+1=\ind(L_2)$. We further classify the unfolding families with  more than one parameter.  We already defined rank-1 unfolding families by \eqref{eq:rank}.  When there are at least two parameters, a family $\{f_\eps\}$ is called a {\em rank-2 unfolding family }if one of the following conditions is satisfied:
\begin{itemize}[nosep]
\item 
$\rank\left.\dfrac{d(\mu(\eps),\omega_1(\eps))}{d\eps}\right|_{\eps=0}=2,
\quad
\dfrac{\omega_1(0)}{2\pi}\notin \left\{0,\dfrac{1}{2}\right\},
\quad
\dfrac{\omega_{2}(\eps)}{2\pi}\not\equiv \dfrac{p}{q} \quad \mbox{for}\quad  3\leq q \leqslant 8,
$\\[5pt] 
\item 
$\rank\left.\dfrac{d(\mu(\eps),\omega_2(\eps))}{d\eps}\right|_{\eps=0}=2,
\quad
\dfrac{\omega_2(0)}{2\pi}\notin \left\{0,\dfrac{1}{2}\right\},
\quad
\dfrac{\omega_{1}(\eps)}{2\pi}\not\equiv \dfrac{p}{q} \quad \mbox{for}\quad  3\leq q \leqslant 8,
$\\[5pt] 
\item 
$
\rank\left.\dfrac{d(\mu(\eps),\theta(\eps))}{d\eps}\right|_{\eps=0}=2,
\quad
\left(\dfrac{\omega_1(0)}{2\pi},\dfrac{\omega_2(0)}{2\pi}\right)\notin \left\{0,\dfrac{1}{2}\right\}^2,
\quad
\dfrac{\omega_{1,2}(\eps)}{2\pi}\not\equiv \dfrac{p}{q} \quad \mbox{for} \quad  3\leq q \leqslant 8.
$
\end{itemize}
Here $p,q$ are coprime integers, and `$\not\equiv$' means not equal in a small neighborhood of $\eps=0$. Recall that we take $\omega_{1,2}\in [0,\pi]$ (see \eqref{eq:argu}). Hence, the multipliers $\lambda$ and $\gamma$ are real if and only if
$\omega_1/2\pi \in \{0,1/2\}$ and $\omega_2/2\pi \in \{0,1/2\}$, respectively. 

If the family has least three parameters, we call it a {\em rank-3 unfolding family} if one of the following conditions is satisfied:
\begin{itemize}[nosep]
\item 
$
\rank\left.\dfrac{d(\mu(\eps),\omega_1(\eps),\theta(\eps))}{d\eps}\right|_{\eps=0}=3,\quad
\dfrac{\omega_1(0)}{2\pi}\notin \left\{0,\dfrac{1}{2}\right\},
\quad
\dfrac{\omega_{2}(\eps)}{2\pi}\not\equiv \dfrac{p}{q} \quad \mbox{for}\quad 3\leq |q| \leqslant 8,
$\\[5pt] 
\item 
$
\rank\left.\dfrac{d(\mu(\eps),\omega_2(\eps),\theta(\eps))}{d\eps}\right|_{\eps=0}=3,\quad
\dfrac{\omega_2(0)}{2\pi}\notin \left\{0,\dfrac{1}{2}\right\},
\quad
\dfrac{\omega_{1}(\eps)}{2\pi}\not\equiv \dfrac{p}{q} \quad \mbox{for}\quad 3\leq |q| \leqslant 8,
$\\[5pt] 
\item   $
\rank\left.\dfrac{d(\mu(\eps),\omega_1(\eps),\omega_2(\eps))}{d\eps}\right|_{\eps=0}=3,\quad
\left(\dfrac{\omega_1(0)}{2\pi},\dfrac{\omega_2(0)}{2\pi}\right)\notin \left\{0,\dfrac{1}{2}\right\}^2.
$
\end{itemize}
Note that the rank-2 and rank-3 unfolding families are generic in the space of families with at least two and, respectively, three parameters.

Recall that two hyperbolic periodic points $P$ and $Q$ are {\em homoclinically related} if they have the same index and their invariant manifolds intersect cyclically and transversely: $W^s(P)\pitchfork W^u(Q)\neq \emptyset$ and $W^u(P)\pitchfork W^s(Q) \neq \emptyset$. Two hyperbolic basic sets are homoclinically related if they have two periodic orbits which are homoclinically related.
\begin{thma}\label{thm:unfold}
Let $f\in \diff ^r(\mathcal{M})$, with $r=3,\dots,\infty,\omega$,  have a non-degenerate heterodimensional cycle such that at least one of $\lambda$ and $\gamma$ is nonreal\footnote{Here we require one more regularity for $f$, because this guarantees a smooth dependence on the parameters of the coordinates used in the proof. For $r=2$, the result is still true with a minor adjustment to the arguments, see the discussion above \citep[Corollary 2]{LT21}  for details.}. Let $\{f_\eps\}$ be either a rank-$i$ ($i=2,3$) unfolding family, or a rank-1 unfolding family with $f$ satisfying the condition in Corollary \ref{cor:limit}. Given any neighborhood $U$ of the cycle, there exist open regions in the $\eps$-space accumulating on $\eps=0$ such that each corresponding diffeomorphism $f_\eps$   has robust heterodimensional dynamics in $U$, which involve
two non-trivial hyperbolic basic  sets $\Lambda_1$ of index $d_1$ and $\Lambda_2$ of index $d_2$. The sets $\Lambda_1$  and $\Lambda_2$ are  homoclinically related to the continuations of $L_1$ and $L_2$, respectively. Moreover, we have that
\begin{itemize}[nosep]
\item if the family is  rank-1  or -2, then at least one of the following holds: $\Lambda_1$ is a cs-blender,  $\Lambda_2$ is a cu-blender, or both; 
\item if the family is rank-3, then $\Lambda_1$ is a cs-blender and $\Lambda_2$ is a cu-blender.
\end{itemize}
\end{thma}

\begin{rem}
This result reduces the number of parameters necessary to obtain blenders by  1 or 2 depending on the situation, in comparison with that required in \citep[Theorem 7]{LT21}. This makes the task to find blenders in specific examples easier.
In particular, Theorem~\ref{thm:unfold} is used in \citep{Tu:25} to  obtain certain pseudohyperbolic attractors containing robust heterodimensional dynamics within a 2-parameter family. With \citep{Li16,LT17,LT20}, it can also be applied to obtain blenders in the attractors in Lorenz-like systems.
\end{rem}

We stress that the conditions for the families can be verified explicitly and can exist in  restricted settings such as perturbations keeping a symmetry. 
The theorem follows from Theorem \ref{thm:unfolding1} in Section \ref{sec:subcases}, which concerns one-parameter unfoldings of heterodimensional cycles with different arithmetic properties of the moduli. It should be noted that not for all cases the cycle needs to be unfolded. For example, in the second case of rank-3 unfolding families, the result is achieved by keeping $\mu=0$ and adjusting $\omega_i$ so that  $\omega_i/2\pi$ are both irrational. If we do not require the homoclinic relations between $\Lambda_i$ and the continuations of $L_i$, then the cycle can be preserved also in the first case of rank-3 unfolding families. This fact can be inferred from Theorem \ref{thm:df}.

\subsection{A dichotomy for local dynamics near heterodimensional cycles}
In closing the introduction, we present a result on the complexity of the dynamics near a heterodimensional cycle. To that aim, we need the bifurcation results for cycles with real $\lambda$ and $\gamma$ obtained  in~\citep{LT21}, which can be summarised as 
\begin{thms}[{\citep[Theorem 2 and Corollary 2]{LT21}}]\label{thm:type1.2}
Let $f$ have a non-degenerate  heterodimensional cycle with both real central multipliers.
Suppose that
\begin{itemize}[nosep]
\item when the cycle is type-I, $\{f_\eps\}$ is rank-2  or rank-1 with irrational $\theta(0)$; and
\item when the cycle is type-II, $\{f_\eps\}$ is rank-2 .
\end{itemize}
Given any neighborhood $U$ of the cycle, there exist open regions in the $\eps$-space accumulating on $\eps=0$ such that each corresponding diffeomorphism $f_\eps$   has robust heterodimensional dynamics in $U$, which involves
two non-trivial hyperbolic basic  sets $\Lambda_1$ of index $d_1$ and $\Lambda_2$ of index $d_2$. 
\end{thms} 


Note from Theorem \ref{thm:rational} that by taking the union of $\mathcal{H}'_f\subset \mathcal{H}_f$ over different rational values of the moduli, we obtain  a dense subset $\mathcal{S}_f\subset \mathcal{H}_f$ where each diffeomorphism $g$ has {\simple} defined in Section \ref{sec:moduli}. Since non-degenerate cycles are open and dense in $\mathcal{H}_f$ (see Appendix~\ref{sec:nd}), the above  observation along with Corollary \ref{cor:limit}, Theorem   \ref{thm:unfold} and Theorem \ref{thm:type1.2} leads to
\begin{cor}\label{cor:dich}
Let $f\in \diff ^r(\mathcal{M})$, with $r=3,\dots,\infty,\omega$, have a heterodimensional cycle $\Gamma$. There exists a dense subset $\mathcal{S}\in \mathcal{H}_f$, and, when the cycle is not type-I, there also exists  a dense subset $\mathcal{C}\in \mathcal{H}_f\setminus \mathcal{S}$, such that
\begin{itemize}[nosep]
\item any close arc transverse to $\mathcal{H}_f$ at $g\in\mathcal{S}$ has {\simple} near the cycle $\Gamma_g$, and
\item for $g\in\mathcal{C}$, the cycle $\Gamma_g$ is a limit of blenders, and  any close arc transverse to $\mathcal{H}_f$ at $g$ intersects a sequence of open sets accumulating on $g$, which correspond to diffeomorphisms having robust heterodimensional dynamics involving blenders.
\end{itemize}
\end{cor}

\subsection{Plan of the proofs}

In Section \ref{sec:subcases}, to facilitate the discussion of the problem, we first define two types of heterodimensional cycles, {\em saddle-focus} and {\em double-focus}, depending on whether only one or both central multipliers are nonreal. Then, according to the types and finer arithmetic properties of $\theta,\omega_1/2\pi,\omega_2/2\pi$ (e.g. rational independence), we split Theorem \ref{thm:main} into  slightly more technical parts, Theorems \ref{thm:sf}--\ref{thm:rational2}. After that, based on these theorem, we state a bifurcation result, Theorem \ref{thm:unfolding1}, which implies Theorem \ref{thm:unfold}.

In Section \ref{sec:firstre}, we introduce the first return maps associated with heterodimensional cycles, and present the cross-form formulas for these maps that are obtained in \citep{LT21} in a more unified way. The  proofs of  our results are mostly done by analyzing these maps.

The main results,   Theorems \ref{thm:sf}--\ref{thm:rational2}, are proved in  Sections \ref{sec:blendersf} and \ref{sec:blenderdf} for saddle-focus, and, respectively, double-focus cycles. A key step is the use of a criterion for the existence of blenders, which is adapted to the cross-form formulas, such that the problem of finding blenders is reduced to finding sequences $\{(k_n,m_n)\}\subset \mathbb{N}^2$ such that the set
$$\left\{ m_n-k_n\theta, k_n \dfrac{\omega_1}{2\pi},m_n \dfrac{\omega_2}{2\pi} \right\}\subset \mathbb{R}\times S^1\times S^1$$
has suitable accumulation points. The criterion is introduced in Section \ref{sec:crit}.

We deal with Theorem \ref{thm:unfolding1} and Theorem \ref{thm:rational} in Sections~\ref{sec:unfolding} and \ref{sec:rational}, respectively. The proofs there follow closely the corresponding proofs in \citep{LT21}, with the use of the newly obtained non-perturbative results Theorems \ref{thm:sf}--\ref{thm:rational2}.


\section{Results for different types of heterodimensional cycles}\label{sec:subcases}
In this section, we state several detailed results  considering more subcases of the arithmetic properties of the moduli. These results  imply Theorems \ref{thm:main} and \ref{thm:unfold}.  We distinguish three main types of heterodimensional cycles, defined as follows. Let $f\in \diff^r(\mathcal{M})$ $(r=2,\dots,\infty)$  have a heterodimensional cycle involving two hyperbolic periodic orbits $L_1$ and $L_2$ satisfying
$\ind(L_1)+1=\ind(L_2).$  Denoting $d=\dim (\mathcal{M})$ and $d_i=\ind(L_i)$ $(i=1,2)$, we can order the multipliers of $L_i$ as
\begin{equation}\label{eq:mult}
|\lambda_{i,d-d_i}|\leqslant\dots\leqslant |\lambda_{i,2}|\leqslant |\lambda_{i,1}|\leqslant 1\leqslant |\gamma_{i,1}|\leqslant |\gamma_{i,2}|\leqslant \dots\leqslant |\gamma_{i,d_i}|.
\end{equation}
Up to an arbitrarily  $C^r$-small perturbation, we can assume that the multipliers are in general position. So, the center-stable multiplier(s)  either  is $\lambda_{1,1}$ being simple and real, or are $\lambda_{1,1}$ and $\lambda_{1,2}=\lambda^*_{1,1}$ being a pair of simple complex conjugate numbers; and the center-unstable multiplier(s)   either  is $\gamma_{2,1}$ being simple and real or are $\gamma_{2,1}$ and $\gamma_{2,2}=\gamma^*_{2,1}$ being a pair of simple complex conjugate numbers. We say that the heterodimensional cycle is of
\begin{itemize}[nosep]
\item saddle type, if $\lambda_{1,1}=:\lambda$ and $\gamma_{2,1}=:\gamma$ are real and satisfy $|\lambda_{1,2}|<|\lambda_{1,1}|$ and 
$|\gamma_{2,1}|<|\gamma_{2,2} |$;
\item saddle-focus type, if either 
$$\lambda_{1,1}=\lambda_{1,2}^*=\lambda e^{i\omega}, \;\omega\in(0,\pi), \quad\mbox{and}\quad \gamma_{2,1}=:\gamma\quad \mbox{is real},$$
where $\lambda>|\lambda_{1,3}|$ and $|\gamma_{2,1}|<|\gamma_{2,2} |$, or
$$\gamma_{2,1}=\gamma_{2,2}^*=\gamma e^{i\omega},\;\omega\in(0,\pi), \quad\mbox{and}\quad \lambda_{1,1}=:\lambda\quad \mbox{is real},$$
where $\gamma<|\gamma_{2,3}|$ and $|\lambda_{1,2}|<|\lambda_{1,1}|$;
\item double-focus type, if
$$\lambda_{1,1}=\lambda_{1,2}^*=\lambda e^{i\omega_1},\;\omega_1\in(0,\pi),\quad\mbox{and}\quad 
\gamma_{2,1}=\gamma_{2,2}^*=\gamma e^{ i\omega_2},\;\omega_2\in(0,\pi),$$
where $\lambda>|\lambda_{1,3}|$ and $0<\gamma<|\gamma_{2,3}|$.
\end{itemize}
As the second case for saddle-focus type can be reduced to
the first one by considering $f^{-1}$, we always assume that for saddle-focus cycles $\lambda_{1,1}$ is nonreal and $\gamma_{2,1}$ is real. 


As mentioned, in this paper we focus on   saddle-focus and double-focus heterodimensional cycles (while the saddle ones have been comprehensively studied in \citep{LT21}). To formulate our results in a compact way, we use a version of the notion of activation introduced in \citep{BD08}.
\begin{defi}[Activation]\label{defi:act}
A cs-blender $\Lambda$ of index $i$ is said to be {\em activated} by a hyperbolic set $H$ if $\ind(H)=i+1$, $W^u(H)$ intersects $ W^s(\Lambda)$ transversely, and $W^s(H)$ contains a disc in the open set $\mathcal{D}$  in Definition \ref{defi:blender} (with $d^{ss}={d-i-1}$) so that the non-transverse intersection $W^s(H)\cap W^u(\Lambda)$ is $C^1$-robustly non-empty.\\
Similarly, a cu-blender $\Lambda$ of index $i$ is {\em activated} by a hyperbolic set $H$ if  $\ind(H)=i-1$,  $W^s(H)\pitchfork W^u(\Lambda)\neq\emptyset$, and  $W^u(H)$ contains a disc in  an open set $\mathcal{D}$ of embeddings of $(i-1)$-dimensional discs, each of which intersects $C^1$-robustly the set $W^s(\Lambda)$. 
\end{defi}

Let us start with saddle-focus heterodimensional cycles. As mentioned, we  always assume that the center-stable multipliers are a pair of complex conjugate numbers and the center-unstable multiplier is real. Throughout the remaining paper, we  use the notation $p/q$ to refer to rational numbers with $p,q$ being coprime integers; and, in all cases, a cs-blender always has index $d_1$ and a cu-blender always has index $d_2$.

\begin{thm}\label{thm:sf}
Let $f\in \diff^r(\mathcal{M})$ have a non-degenerate saddle-focus heterodimensional cycle, and let $U$ be any small neighborhood of the cycle. 
\begin{enumerate}[label=(\arabic*),nosep]
\item If $\omega/2\pi$ is irrational, then there exists a cu-blender in $U$, homoclinically related to $L_2$.
\item  If $\omega/2\pi$ is irrational, and $\theta,\omega/2\pi,1$ are rationally independent (so $\theta$ is also irrational), then there exists a cs-blender in $U$, in addition to the cu-blender of case (1) which is homoclinically related to $L_2$. The two blenders activate each other, and in particular $L_2$ activates the cs-blender.
\item If   $\omega/2\pi=p/q$ with $q\geqslant 9$ and $\theta$ is irrational, then there  exists a blender in $U$.
\end{enumerate}
\end{thm}

\begin{rem}
The type of the blender  in case (3) depends on certain coefficients of the first return map, see case 1 of Proposition \ref{prop:blendersf}.
\end{rem}

The theorem is proved in Section \ref{sec:blendersf}. Case (2) here, as well as case (2) in Theorem \ref{thm:df} below, follows from \citep[Theorem 7]{LT21}, but we reprove them with a more general approach. Let us explain briefly how the lower bound for $q$ appears in case (3). The existence of blenders is determined by the central dynamics of the first return maps $T_{k,m}$ along a heterodimensional cycle, where $k$ and $m$ represent the number of iterations spent near the two periodic orbits $L_1$ and $L_2$, respectively. When the center-stable multiplier of $L_1$ is nonreal, the rotation brought by it leads to the emergence of trigonometric terms  like $\sin(k\omega+\eta)$ in the formula for the one-dimensional central part of the first return map. To prove the existence of blenders, we need  certain inequalities involving these trigonometric terms to be satisfied (see \eqref{eq:bad}), and in most cases we only need $\sin(k\omega+\eta)$ to take sufficiently many different values for different $k$. This is ensured by a lower bound on $q$. See Section \ref{sec:fr_sf} for details.

For double-focus heterodimensional cycles, we have 
\begin{thm}\label{thm:df}
Let $f\in \diff^r(\mathcal{M})$ have a non-degenerate double-focus heterodimensional cycle, and let $U$ be any small neighborhood of the cycle.
\begin{enumerate}[label=(\arabic*),nosep]
\item Suppose that at least one of $\omega_i/2\pi$ $(i=1,2)$ is irrational. 
\begin{enumerate}[label=(\roman*),nosep]
\item If $\omega_1/2\pi$ is irrational and $\omega_2/2\pi=p/q$ with $q\geqslant \q2$, then there exists  a cu-blender in $U$, homoclinically related to $L_2$.
\item If $\omega_2/2\pi$ is irrational and $\omega_1/2\pi=p/q$ with $q\geqslant \q2$, then there exists  a cs-blender in $U$, homoclinically related to $L_1$.
\item If both $\omega_1/2\pi$ and $\omega_2/2\pi$ are irrational,  then in $U$ there exist simultaneously a cs-blender, homoclinically related to $L_1$, and a cu-blender, homoclinically related to $L_2$; the two blenders activate each other, and in particular $L_1$ activates the cu-blender and $L_2$ activates the cs-blender.
\end{enumerate}

\item  Suppose that $\theta$ and exactly one of $\omega_i/2\pi$  are irrational.
\begin{enumerate}[label=(\roman*),nosep]
\item If $\theta,\omega_1/2\pi,1$ are rationally independent and $\omega_2/2\pi=p/q$ with $q\geqslant \q2$, then   there exists in $U$ a cs-blender, in addition to the cu-blender of (1.i) which is homoclinically related to $L_2$.
\item If $\theta^{-1},\omega_2/2\pi,1$ are rationally independent and $\omega_1/2\pi=p/q$ with $q\geqslant \q2$, then  there exists a cu-blender in $U$, in addition to the cs-blender of (1.ii) which is  homoclinically related to $L_1$.
\end{enumerate}
In each of the two cases, the two blenders activate each other, and in particular $L_2$ activates the cs-blender in (2.i) and $L_1$ activates the cu-blender in (2.ii).

\item If $\theta$ is irrational, and $\omega_1/2\pi=p_1/q_1$ and $\omega_2/2\pi=p_2/q_2$ such that either $q_1\geqslant 9$ and $q_2\geqslant \q2$ or $q_2\geqslant 9$ and $q_1\geqslant \q2$, then there always exists a blender in $U$.
\end{enumerate}
In all cases, a cs-blender always has index $d_1$ and a cu-blender always has index $d_2$.
\end{thm}
\begin{rem}
The type of the blender  in case (3) depends on certain coefficients of the first return map, see case 1 of Proposition \ref{prop:blenderdf}.
\end{rem}

This theorem is proved in Section \ref{sec:blenderdf}, which along with   Theorem \ref{thm:sf}  implies
Theorem \ref{thm:main}. The reason for the lower bounds of $q$ is the same as for the saddle-focus cycles, where the related inequalities are given in \eqref{eq:bad2}.

The next result concerns the possibilities for the simultaneous appearance of a cs-blender and a cu-blender in cases (3) of Theorems \ref{thm:sf} and \ref{thm:df}. In fact, such coexistence is ensured when $q$ is larger than  a  number which depends on the diffeomorphism and is locally constant.   Recall that, for a diffeomorphism $g$ having a non-degenerate cycle $\Gamma$, we denote by $\mathcal{H}_{g}$ the codimension-1 subspace in $\diff^r(\mathcal{M})$ where every $f\in \mathcal{H}_{g}$ has a heterodimensional cycle $\Gamma_f$ near $\Gamma$. Also note that being non-degenerate is $C^r$-open in $\mathcal{H}_{g}$,  see Appendix~\ref{sec:nd}.

\begin{thm}\label{thm:rational2}
Let $g$ have a non-degenerate heterodimensional cycle of saddle-focus or double-focus type $\Gamma$. For any sufficiently small neighborhood  $\mathcal{U}$ of $g$ in $\mathcal{H}_{g}$, there exists a number $N(\mathcal{U})$ such that a diffeomorphism $f\in\mathcal{U}$  has a pair of mutually activating cs-blender and cu-blender in any neighborhood of $\Gamma_f$ if
\begin{itemize}[nosep]
\item  in the saddle-focus case, $\theta(f)$ is irrational and $\omega(f)/2\pi=p/q$ with $q \geq N(\mathcal{U})$, or
\item  in the double-focus case, $\theta(f)$ is irrational and $\omega_i(f)/2\pi=p_i/q_i$ with either $q_1 \geq N(\mathcal{U})$ and $q_2 \geq \q2$ or $q_2 \geq N(\mathcal{U})$ and $q_1 \geq \q2$.
\end{itemize}
\end{thm}  

We prove this theorem in Sections \ref{sec:proofsf} and \ref{sec:proofdf}. In light of this result, by taking $q\to\infty$, it is intuitive to expect coexisting blenders near saddle-focus cycles when both $\theta$ and $\omega/2\pi$ are irrational. It is true if $\theta$, $\omega/2\pi$ and 1 are rationally independent, as in case (1) of Theorem \ref{thm:sf}. However, in the absence of such rational independence, we cannot prove coexisting blenders, see case 2 of Proposition~\ref{prop:blendersf}.  Similar comments  apply to  double-focus cycles, see case 2.2 of Proposition~\ref{prop:blenderdf}. 

%

Let us now present some perturbative results based on the above theorems. The results show that, by splitting the heterodimensional cycle, one builds homoclinic relations between the involved periodic orbits and the blenders found above, and hence obtains robust heterodimensional dynamics.  Recall that a one-parameter family $\{f_\eps\}$ is of rank-1 if $d\mu/d\eps\neq 0$ at $\eps=0$, see \eqref{eq:rank}. Let us denote the cs-blenders and cu-blenders in Theorems \ref{thm:sf}--\ref{thm:rational2} by $\Lambda^{cs}$ and $\Lambda^{cu}$, respectively. We use the same notations $L_{1,2},\Lambda^{cs},\Lambda^{cu}$ for their continuations at $\eps\neq 0$.

\begin{thm}\label{thm:unfolding1}
Consider any   rank-1 unfolding family  $\{f_\eps\}$ of $C^r$ ($r=3,\dots,\infty,\omega$) diffeomorphisms with  $f_0=f$ from either of Theorems \ref{thm:sf}--\ref{thm:rational2}. Let $U$ be any neighborhood of the cycle of $f$. There exists a sequence $\{\eps_j\}$ accumulating on $\eps=0$ such that for each $f_{\eps_j}$,

\begin{itemize}[nosep]
\item in Theorem \ref{thm:sf}.(1) and Theorem \ref{thm:df}.(1.i), $\Lambda^{cu}$ is homoclinically related to $L_2$, and is activated by a non-trivial hyperbolic basic set in $U$ which is homoclinically related to  $L_1$;
\item in  Theorem \ref{thm:df}.(1.ii), $\Lambda^{cs}$ is homoclinically related to $L_1$, and is activated by a non-trivial hyperbolic basic set  in $U$ which is homoclinically related to  $L_2$;
\item in Theorem \ref{thm:sf}.(2), Theorem \ref{thm:df}.(2) and Theorem \ref{thm:rational2}, $\Lambda^{cs}$ and $\Lambda^{cu}$ activate each other, and are homoclinically related to $L_1$ and $L_2$, respectively;
\item in Theorem \ref{thm:sf}.(3) and Theorem \ref{thm:df}.(3), if the blender is center-stable, it is homoclinically related to $L_1$, and is activated by a non-trivial hyperbolic basic set  in $U$ which is homoclinically related to  $L_2$; if  the blender is center-unstable, then it is homoclinically related to $L_2$, and is activated by a non-trivial hyperbolic basic set in $U$ which is homoclinically related to  $L_1$.
\end{itemize}
Moreover, in all cases, the 
periodic orbits, blenders and hyperbolic sets involved in the robust heterodimensional dynamics, as well as the heteroclinic orbits constituting the intersections involved in the  activations and homoclinic relations all lie in $U$.
\end{thm}

This theorem is proved in Section \ref{sec:unfolding}. Theorem \ref{thm:unfold} follows immediately from the above results, by using the additional parameters to change $\theta,\omega,\omega_{1,2}$ to fulfill the hypotheses required in corresponding cases. The last part of  Theorem \ref{thm:unfolding1} shows that all the objects related to the obtained robust heterodimensional dynamics are localized arbitrarily close to the original cycle. This fact is important to the (local) stabilization theory of heterodimensional cycles \citep{BDK12,LT21,BDL}.


\section{First return maps}\label{sec:firstre}
In this section, we present the formulas obtained in \cite{LT21} for the first return maps of non-degenerate saddle-focus and double-focus heterodimensional cycles. The proofs of our main results are based on the analysis of these formulas.

\subsection{Local map and transition map}\label{sec:setting}
Recall that $f\in \diff^r(\mathcal{M})$ with $r=2,\dots,\infty,\omega$ has a heterodimensional cycle involving two hyperbolic periodic orbits $L_1$ and $L_2$, where their indices $d_1$ and $d_2$ satisfy $d_1+1=d_2$. Let $\Gamma^0$ be a heteroclinic orbit in $W^u(L_1)\cap W^s(L_2)$, and $\Gamma^1$ be a heteroclinic orbit in $W^s(L_1)\cap W^u(L_2)$. We further denote the heterodimensional cycle as the closed invariant set
$
\Gamma=L_1\cup L_2\cup \Gamma^0\cup \Gamma^1.
$

Take two points $O_1$ and $O_2$ from the two periodic orbits $L_1$ and $L_2$. For $i=1,2$, let $U_{0i}$ be some sufficiently small neighborhoods  of $O_i$.
We define the {\em local maps} $F_i:=f^{per(O_i)}|_{U_{0i}}$, where $per(O_i)$ are the periods of $O_i$.  
Next, using the heteroclinic orbits $\Gamma^0$ and $\Gamma^1$, we define {\em transition maps} connecting neighborhoods of $O_1$ and $O_2$. Take in $\Gamma^0$ a pair of points $M^-_1\in W^u_{loc}(O_1)$ and $M^+_2\in W^s_{loc}(O_2)$ so that $f^{n_{1}}(M^-_1)=M^+_2$ for some integer $n_1>0$, and in $\Gamma^1$ a pair of points $M^-_2\in W^u_{loc}(O_2)$ and $M^+_1\in W^s_{loc}(O_1)$ so that $f^{n_{2}}(M^-_2)=M^+_1$ for some integer $n_2>0$. The transition maps are then defined as $F_{12}:=f^{n_1}$, taking a small neighborhood of $M^-_1$ to a small neighborhood of $M^+_2$,  and $F_{21}:=f^{n_2}$, taking a small neighborhood of $M^-_2$ to a small neighborhood of $M^+_1$. The first return maps will be defined as compositions of the local map and  transition map.
\begin{rem}\label{rem:mu0}
The splitting parameter introduced in Section~\ref{sec:moduli} can be defined to measure the signed distance between $F_{12}(W^u_{loc}(O_1))$ and $W^s_{loc}(O_2)$ in a small neighborhood of $M^+_2$.
\end{rem}

As mentioned, we assume that the cycle of $f$ is {\em non-degenerate} in the sense of \citep{LT21}, which can always be achieved by an arbitrarily $C^r$-small perturbation. There are four non-degeneracy conditions, whose precise definitions are given in  Appendix~\ref{sec:nd}. Roughly speaking, the first three require that
\begin{itemize}[nosep]
\item $T_1(W^u_{loc}(O_1))$ intersects   transversely the  local extended stable manifold $W^{sE}(O_2)$ at $M^+_2$, and  $T^{-1}_1(W^s_{loc}(O_2))$ intersects transversely the local extended unstable manifold $W^{uE}(O_1)$ at $M^-_1$,
\item $T_2^{-1}(W^s_{loc}(O_1))$ is not tangent to the strong-unstable leaf $\ell^{uu}\subset W^u(O_2)$ at $M^-_2$, and  $T_2(W^u_{loc}(O_2))$ is not tangent to the strong-stable leaf $\ell^{ss}\subset W^s(O_1)$ at $M^+_1$, and
\item $M^+_1$ does not belong to the strong-stable manifold $W^{ss}_{loc}(O_1)$ and $M^-_2$ does not belong to the strong-unstable manifold $W^{uu}_{loc}(O_2)$.
\end{itemize}

\begin{figure}[!h]
\begin{center}
\includegraphics[width=0.8\textwidth]{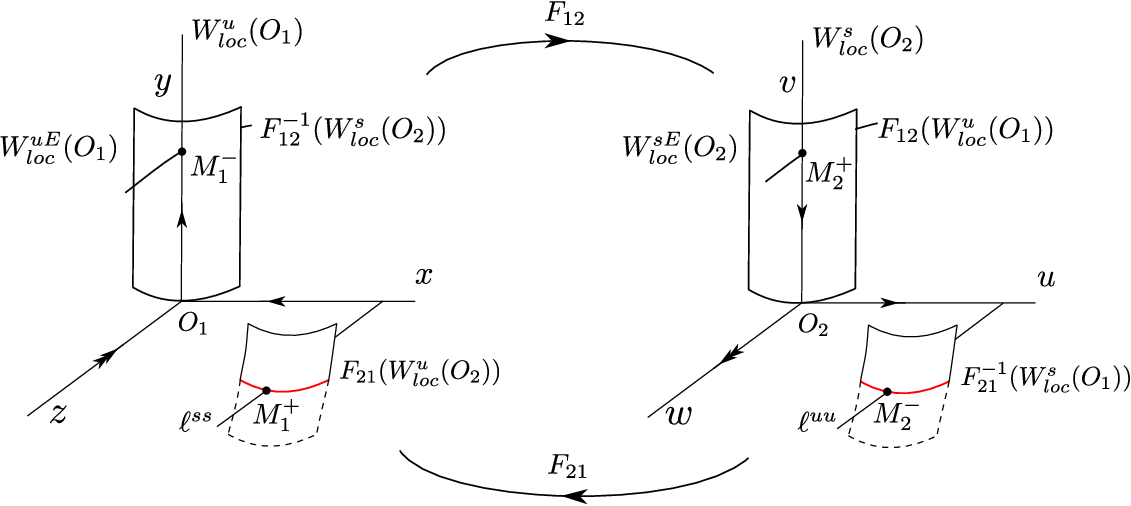}
\end{center}
\caption{Here $W^u_{loc}(O_1)$ is the  $y$-axis, $W^s_{loc}(O_1)$ is the  $(x,z)$-plane, $W^s_{loc}(O_2)$ is the  $v$-axis, and $W^u_{loc}(O_2)$ is the  $(u,w)$-plane. The two red curves are $\ell_1=W^s_{loc}(O_1)\cap F_{21}(W^u_{loc}(O_2))$ and $\ell_2=F^{-1}_{21}(\ell_1)$.}
\label{fig:gc}
\end{figure}

Here  $\ell^{uu},\ell^{ss},W^{ss}_{loc}(O_1), W^{uu}_{loc}(O_2)$ are  defined in the usual way; $W^{sE}_{loc}(O_2)$ is the invariant manifold corresponding to the stable and weak-unstable directions at $O_2$, and similarly for $W^{uE}_{loc}(O_1)$. Figure~\ref{fig:gc} provides an illustration.
  The last non-degeneracy condition is more delicate, requiring that the projections of the transverse intersection points $M^+_1$ and $M^-_2$ to the central directions do not lie in certain lines.  See Appendix~\ref{sec:nd} for details.

\subsection{Saddle-focus case}\label{sec:fr_sf}
The first return map of a saddle-focus heterodimensional cycle is defined as the composition  $T_{k,m}:=F_{21}\circ F^m_2\circ F_{12}\circ F^k_1$, which, for suitable choices of $(k,m)$, takes points in a small neighborhood of $M^+_1$ back to the neighborhood.  It is shown in \citep[Sections 6.1.1--6.1.3]{LT21} that there exist coordinates
 $(X,Y,Z)\in \mathbb{R}\times\mathbb{R}^{d_1}\times\mathbb{R}^{d-d_1-1}$ near  $M^+_1$ such that the first return map can be expressed in a cross-form, which manifests the partial hyperbolicity, for a collection of pairs $(k,m)$. In what follows we introduce this formula. 

Recall that we  denote the pair of complex conjugate center-stable multipliers by $\lambda e^{\pm i\omega}$ $(0<\omega<\pi)$ and the real center-unstable multiplier by $\gamma$. 
For $\delta>0$, define 
\begin{equation}\label{eq:complex:domain}
\Pi=[-\delta,\delta]\times[-\delta,\delta]^{d_1}\times [-\delta,\delta]^{d-d_1-1},
\end{equation}
where $d_1=\ind(O_1)$. 

It is proven  that
 there exist coefficients $A,B,b,\eta_{1,2}\in\mathbb{R}$ (which are expressed in terms of the derivatives $\D F_{12}|_{M^-_1}$ and $\D F_{21}|_{M^+_2}$), $u^-\in\mathbb{R}$ (which is given by the position of $M^-_2$, see \eqref{eq:intro:2}) and a positive number $C(\delta)$ depending  on $\delta$ such that, for all sufficiently small $\delta$ and sufficiently large $k$ and $m$ satisfying
\begin{equation}\label{eq:crosscon}
|\lambda^k\gamma^m B\sin (k\omega +\eta_2)-bu^-|<C(\delta),
\end{equation}
one has that $(\bar X,\bar Y,\bar Z):=T_{k,m}(X,Y,Z)\in\Pi$ for a point $(X,Y,Z)\in \Pi$ 
if and only if
\footnote{Equation~\eqref{eq:complex:T_k,m_cross_0mu} is obtained from   \citep[Equation  (6.18)]{LT21} by dropping the subscript of $X_1$ and grouping the coordinates $(X_2,Z)=:Z$.}
\begin{equation}\label{eq:complex:T_k,m_cross_0mu}
\begin{aligned}
\bar X &=  A_{k,m} X + B_{k,m} +\phi_{1}(X,\bar Y,Z),\\
 Y &=\phi_2(X,\bar Y,Z),\quad \bar Z = \phi_3(X, \bar Y,Z),
\end{aligned}
\end{equation}
where the functions $\phi$ satisfy 
$$\phi_1=o(\delta)+o(1)_{k,m\to\infty},\quad
\dfrac{\partial \phi_1}{\partial (X,\bar Y,Z)}=o(1)_{\delta\to 0}+o(1)_{k,m\to\infty},\quad
\|\phi_{2,3}\|_{C^1}=o(1)_{k,m\to\infty},$$
and
\begin{equation}\label{akmbkm}
A_{k,m}=\lambda^k\gamma^m  A \sin(k\omega + \eta_1),\qquad
B_{k,m}=\lambda^k\gamma^m B\sin(k\omega + \eta_2)- b u^-.
\end{equation}
 The non-degeneracy conditions ensure (see \citep[Section 6.1.3]{LT21}) that $A> 0,B> 0,b\neq 0,u^-\neq 0$ and 
\begin{equation}\label{eq:eta}
\tan\eta_1\neq\tan\eta_2.
\end{equation}

The partial hyperbolicity of $T_{k,m}$ follows immediately from a computation of the invariant cone fields using \eqref{eq:complex:T_k,m_cross_0mu}. Denote by $(\Delta X,\Delta Y,\Delta Z)$ the vectors in the tangent spaces.
\begin{lem}[{\citep[Lemma 6.1]{LT21}}]\label{lem:conefields} Given any $K>0$, one can choose $\delta$ sufficiently small such that, for all sufficiently large $k,m$ for which
$\sin(k\omega+\eta_1)$ and $\sin(k\omega+\eta_2)$ are bounded away from zero, the cone fields 
\begin{align}
&\mathcal{C}^{cu} = \{(\Delta X,\Delta Y,\Delta Z): \|\Delta Z\| \leq K (|\Delta X|+\|\Delta Y\|)\},\nonumber\\
&\mathcal{C}^{uu} = \{(\Delta X,\Delta Y,\Delta Z): \max\{|\Delta X|,\|\Delta Z\|\}\leq K  \|\Delta Y\|\},\label{eq:uucone}
\end{align}
which are defined on $\Pi$, are forward-invariant in the sense that if a point $M\in \Pi$ has its image $\bar M=T_{k,m}(M)$ in $\Pi$, then the cone at $M$ 
is mapped strictly into the cone at $\bar M$ by $\D T_{k,m}$; the cone fields
\begin{align}
&\mathcal{C}^{cs} = \{(\Delta X,\Delta Y,\Delta Z): \|\Delta Y\|\leq K (|\Delta X|+\|\Delta Z\|)\},\nonumber\\
&\mathcal{C}^{ss}=\{(\Delta X,\Delta Y,\Delta Z):\max\{|\Delta X|,\|\Delta Y\|\}\leq K  \|\Delta Z\|\},\label{eq:sscone}
\end{align}
which are defined on $\Pi$, are backward-invariant in the sense that if a point $\bar{M}\in\Pi$ has its preimage $M=T^{-1}_{k,m}(\bar M)$ in $\Pi$, 
then the cone at $\bar{M}$ is mapped strictly into the cone at $M$ by $\D T^{-1}_{k,m}$. Moreover, vectors in $\mathcal{C}^{uu}$ and, if $|A_{k,m}|>1$, also vectors in $\mathcal{C}^{cu}$ are expanded by $\D T_{k,m}$; vectors in $\mathcal{C}^{ss}$ and, if $|A_{k,m}|<1$, also vectors in $\mathcal{C}^{cs}$ are contracted by $\D T_{k,m}$.
\end{lem}

\subsection{Double-focus case}\label{sec:fr_df}
In this case, we consider the first return map $T_{k,m}:=F^m_2\circ F_{12}\circ F^k_1\circ F_{21}$ near $M^-_2$. A cross-form formula for  $T_{k,m}$ is obtained in \citep[Sections 6.4.1--6.4.3]{LT21} with coordinates
  $(U,V,W)\in \mathbb{R}\times\mathbb{R}^{d-d_2}\times\mathbb{R}^{d_2-1}$ in a small neighborhood  $M^-_2$, which we describe below. For simplicity, many coefficients here will carry the same notations as those in the saddle-focus case, but they are different coefficients. This will not cause ambiguity as we never discuss  the two cases at the same time in the proofs.

Recall that, for a double-focus heterodimensional cycle, we denote the center-stable multipliers by $\lambda e^{\pm i\omega_1}$ $(0<\omega_1<\pi)$ and the center-unstable multipliers by $\gamma e^{\pm i\omega_2}$ $(0<\omega_2<\pi)$.  For $\delta>0$, define
\begin{equation}\label{eq:df:domain}
\Pi=[-\delta,\delta]\times[-\delta,\delta]^{d-d_2}\times [-\delta,\delta]^{d_2-1},
\end{equation}
where $d_2=\ind(O_2)$. 

It is proven that there exist coefficients $A^*,B^*,b,\eta_{1,2,3}$ (which are expressed in terms of the derivatives $\D F_{12}|_{M^-_1}$ and $\D F_{21}|_{M^+_2}$), $u^-_{1,2}$ (which are given by the position of $M^-_2$, see \eqref{eq:intro:2}) and a positive number $c(\delta)=o(1)_{\delta\to 0}$ depending  on $\delta$ such that, for all sufficiently small $\delta$ and sufficiently large $k$ and $m$ satisfying
\begin{equation}\label{eq:quant}
\min\{|\cos m\omega_2|,\,
|\cos(m\omega_2+\eta_3)|,\,
|\cos(m\omega_2+\eta_3) + b\sin(m\omega_2+\eta_3)|\} >c(\delta),
\end{equation}
one has that   $(\bar U,\bar V,\bar W):=T_{k,m}(U,V,W)\in\Pi$ for a point $(U,V,W)\in \Pi$ if and only if
\footnote{Equation~\eqref{eq:df:T_km_new} is obtained from  \citep[Equation (6.44)]{LT21} by dropping the subscript of $U_1$, grouping the coordinates $(U_2,W)=:W$, and dividing $(1+b_{41}\tan(m\omega+\eta_3))$ to both sides of (6.44).}
\begin{equation}\label{eq:df:T_km_new}
\begin{array}{l}
 \bar U 
=  A_{k,m}U+B_{k,m} + \phi_1(U,V,\bar W),\\
\bar{V}=\phi_2(U,V,\bar W),\qquad
W=\phi_3(U,V,\bar W),
\end{array}
\end{equation}
where the functions $\phi$ satisfy 
$$\phi_1=o(\delta)+o(1)_{k,m\to\infty},\quad
\dfrac{\partial \phi_1}{\partial (U, V,\bar W)}=o(1)_{\delta\to 0}+o(1)_{k,m\to\infty},\quad
\|\phi_{2,3}\|_{C^1}=o(\delta)+o(1)_{k,m\to\infty},$$
and
\begin{equation}\label{eq:df:AkmBkm}
\begin{aligned}
&A_{k,m}=\lambda^k\gamma^mA(m\omega_2)\sin(k\omega_1+\eta_1),\quad
B_{k,m}=\lambda^k\gamma^m B(m\omega_2)\sin(k\omega_1+\eta_2)
-\xi(m\omega_2)
\end{aligned}
\end{equation}
with 
\begin{equation}\label{eq:short1}
\begin{aligned}
&A(m\omega_2)=\dfrac{ A^*}{\cos(m\omega_2+\eta_3) + b\sin(m\omega_2+\eta_3)},\quad
B(m\omega_2)=\dfrac{ B^*}{\cos(m\omega_2+\eta_3) + b\sin(m\omega_2+\eta_3)},\\[10pt]
&\xi(m\omega_2)=\dfrac{\cos(m\omega_2+\eta_3)u^-_1+\sin(m\omega_2+\eta_3)u^-_2}{\cos(m\omega_2+\eta_3) + b\sin(m\omega_2+\eta_3)}.
\end{aligned}
\end{equation}
 The non-degeneracy conditions ensure that $A^*> 0,B^*> 0,b\neq 0,u^-\neq 0$ and $\tan\eta_1\neq\tan\eta_2$, see \citep[Section 6.4.3]{LT21}.

One sees that formula \eqref{eq:df:T_km_new} assumes the same form as \eqref{eq:complex:T_k,m_cross_0mu}, with replacing $U,V,W$ by $X,Z,Y$, respectively. By a completely parallel computation, one finds the following counterpart to Lemma \ref{lem:conefields}, which characterizes the partial hyperbolicity of $T_{k,m}$.
\begin{lem}\label{lem:conefields_df} 
Given any $K>0$, one can choose $\delta$ sufficiently small, such that for all sufficiently large $k,m$ satisfying \eqref{eq:quant} and \eqref{eq:bad2}, the cone fields
\begin{align}
&\mathcal{C}^{cu} = \{(\Delta U,\Delta V,\Delta W): \|\Delta V\| \leq K (|\Delta U|+\|\Delta W\|)\},\nonumber\\
&\mathcal{C}^{uu} = \{(\Delta U,\Delta V,\Delta W): \max\{|\Delta U|,\|\Delta V\|\}\leq K  \|\Delta W )\|\},\label{eq:uucone_df}
\end{align}
which are defined  on $\Pi$, are forward-invariant in the sense that if a point $M\in \Pi$ has its image $\bar M=T_{k,m}(M)$ in $\Pi$, then the cone at $M$ 
is mapped into the cone at $\bar M$ by $\D T_{k,m}$; the cone fields
\begin{align}
&\mathcal{C}^{cs} = \{(\Delta U,\Delta V,\Delta W): \|\Delta W\|\leq K (|\Delta U|+\|\Delta V\|)\},\nonumber\\
&\mathcal{C}^{ss}=\{(\Delta U,\Delta V,\Delta W):\max\{|\Delta U|,\|\Delta W\|\}\leq K  \|\Delta V\|\},\label{eq:sscone_df}
\end{align}
which are defined  on $\Pi$, are backward-invariant in the sense that if a point $\bar{M}\in\Pi$ has its preimage $M=T^{-1}_{k,m}(\bar M)$ in $\Pi$, 
then the cone at $\bar{M}$ is mapped into the cone at $M$ by $\D T^{-1}_{k,m}$. Moreover, vectors in $\mathcal{C}^{uu}$ and, if $|A_{k,m}|>1$, also vectors in $\mathcal{C}^{cu}$ are expanded by $\D T_{k,m}$; vectors in $\mathcal{C}^{ss}$ and, if $|A_{k,m}|<1$, also vectors in $\mathcal{C}^{cs}$ are contracted by $\D T_{k,m}$.
\end{lem}

\begin{rem}\label{rem:flow}
The results in this paper are solely based on the analysis of the first return maps, so they can be transferred easily to the continuous-time case, i.e., when $f$ is a flow on any manifold of dimension at least four.  Indeed, one just needs to define  $F_i$ for each $i=1,2$ as the Poincar\'e map of a local transverse cross-section at some point $O_i$ of the periodic orbit, and define $F_{12}$ by the trajectories of the flow from a neighborhood of $M^-_1$ to a neighborhood of $M^+_2$, and similarly for $F_{21}$.
\end{rem}


\section{Blenders near saddle-focus heterodimensional cycles}\label{sec:blendersf}
In this section, we prove Theorem \ref{thm:sf} and the first case in Theorem \ref{thm:rational2}. Both theorems can be divided into two parts, where one concerns the existence of blenders and the other concerns the homoclinic relations. We will obtain results for these two parts separately, and then prove the theorems. For the blender part, we use a criterion obtained in  \citep{LT21}.

\subsection{A criterion for the existence of blenders}\label{sec:crit}
Let us first introduce the notion of {\em activating pair} to further characterize the embedding set $\mathcal{D}$ in Definition \ref{defi:blender}. Take a cube 
$$Q=\{(x_1,\dots,x_d)\mid x_i\in I_i\}\subset \mathbb{R}^d,$$
where $I_i$ are closed intervals in $\mathbb{R}$. A $k$-dimensional disc $S$ is said to {\em cross} $Q$ if the intersection $S\cap Q$ is given by $(x_{i_{k+1}},\dots,x_{i_d})=s(x_{i_1},\dots,x_{i_k})$, where  $(i_1,\dots,i_d)$ is some permutation of $(1,\dots,d)$, and $s$ is a smooth function defined on $I_{i_1}\times \dots\times I_{i_k}$. It crosses $Q$ {\em properly} with respect to a cone field $\mathcal{C}$ defined on $Q$ if the tangent spaces of $S\cap Q$ lie in $\mathcal{C}$, i.e., the tangent spaces of $S$ lie in $\mathcal{C}$.

For an invariant cone field $\mathcal{C}$, we use $\dim (\mathcal{C})$ to denote the largest dimension of linear subspaces of the tangent spaces contained in $\mathcal{C}$.
\begin{defi}[Activating pair]\label{defi:actpair}
A pair $(Q ,\mathcal{C})$ consisting of a cube and a cone field is called an {\em activating pair for a cs-blender} $\Lambda$ if $\dim (\mathcal{C})=\dim W^s(\Lambda)-1$ and if any surface $S$ crossing $Q$ properly with respect to $\mathcal{C}$ intersects $W^u(\Lambda)$. The pair $(Q ,\mathcal{C})$ is an {\em activating pair for a cu-blender} $\Lambda$ if $\dim (\mathcal{C})=\dim W^u(\Lambda)-1$ and any surface $S$ crossing $Q$ properly with respect to $\mathcal{C}$ intersects $W^s(\Lambda)$. The cube $Q$ is called an activating domain.
\end{defi}
For example, if $(\Pi,\mathcal{C}^{ss})$, with $\Pi$ in \eqref{eq:complex:domain} and $\mathcal{C}^{ss}$ in \eqref{eq:sscone}, is an activating pair of some cs-blender, then all curves crossing $\Pi$ properly with respect to $\mathcal{C}^{ss}$ belong to the embedding set $\mathcal{D}$ in Definition \ref{defi:blender}. Each such curve is  the graph of some smooth function $(X,Y)=s(Z)$, which is defined on $[-\delta,\delta]^{d-d_1-1}$ and whose tangent spaces lie in $\mathcal{C}^{ss}$.

Let $\delta>0$ and take a cube
$$Q=[-\delta,\delta]\times [-\delta,\delta]^{d_Y}\times [-\delta,\delta]^{d_Z}$$
with $1+d_Y+d_Z=d$. Let $\{T_{n}\}_{n\in \mathbb N}$ be a sequence of maps defined on pairwise disjoint subsets $\{\sigma_n\}$ of $Q$ such that $T_i(\sigma_i)\cap T_j(\sigma_j)\cap Q =\emptyset$ for $i\neq j$. Define $T:\bigcup_n \sigma_n \to \mathbb{R}^d$ by $T(X,Y,Z)=T_n(X,Y,Z)$ if $(X,Y,Z)\in \sigma_n$. Suppose  for every $n\in \mathbb N$ and  $(X,Y,Z)\in Q$, one has that $(\bar X,\bar Y,\bar Z)=T_n (X,Y,Z)$ if and only if 
\begin{equation}\label{eq:blender:1}
\begin{aligned}
\bar X &= A_{n} X + B_{n} + \phi_1(X,\bar{Y},Z;n), \\
Y &=\phi_2(X,\bar{Y},Z;n), \\
\bar Z &= \phi_3(X,\bar{Y},Z;n),
\end{aligned}
\end{equation}
where the constants $A_n$ and $B_n$, and the $C^1$ functions $\phi$  defined on $Q$ are allowed to depend on $\delta$, and we have the estimates
$$
\phi_1=o(\delta)+o(1)_{n\to\infty},\quad
\dfrac{\partial \phi_1}{\partial (X,\bar Y,Z)}=o(1)_{\delta\to 0}+o(1)_{n\to\infty},\quad
\|\phi_{2,3}\|_{C^1}=o(1)_{\delta\to 0}+o(1)_{n\to\infty}.
$$

\begin{prop}[{\citep[Proposition 4.12]{LT21}}]\label{prop:blender}
Suppose there exists $\delta_0>0$ such that the constants $A_n$ and $B_n$ satisfy that
\begin{itemize}[nosep]
\item[1.] for each $\delta\in(0,\delta_0)$, the set $\{B_n\}$ is dense in an interval $I$ which contains zero and is independent of $\delta$;
\item[2.] for all $\delta\in(0,\delta_0)$ and all $n\in\mathbb N$,
\begin{itemize}[nosep]
\item[2.1.]either $0<C_1<|A_n|<C_2<1$ with $C_{1,2}$ independent of $n$,
\item[2.2.]or $1<C_3<|A_n|<C_4<\infty$ with $C_{3,4}$ independent of $n$.
\end{itemize}
\end{itemize}
Then, there exist $c\in(0,1)$ and 
\begin{equation}\label{eq:Q}
Q'=[-c\delta,c\delta]\times [-\delta,\delta]^{d_Y}\times [-\delta,\delta]^{d_Z},
\end{equation}
such that for every $\delta\in(0,\delta_0)$ one can find a finite set $\mathcal{N}\subset \mathbb{N}$ with arbitrarily large elements such that  the restriction $T_{\mathcal{N}}$ of $T$ to $\bigcup_{n\in\mathcal{N}} \sigma_n\cap Q$ has a blender, which consists of points whose orbits never leave $Q$ under $T_{\mathcal{N}}$, and has  an activating pair $(Q',\mathcal{C})$ with the cone field $\mathcal{C}$  independent of $\delta$ and $n$. Moreover,
\begin{itemize}[nosep]
\item if condition 2.1 holds, $\Lambda$ is a cs-blender of index $d_Y$ and cones in $\mathcal{C}$ are around $Z$-coordinates;
\item if condition 2.2 holds, $\Lambda$ is a cu-blender of index $(d_Y+1)$ and cones in $\mathcal{C}$ are around $Y$-coordinates.
\end{itemize}
\end{prop}

In our case, the maps $T_n$ will be the return maps $T_{k_n,m_n}$ for some sequence of pairs satisfying  conditions 1 and 2.  For completeness we sketch the proof of this proposition in Appendix~\ref{sec:blendercon}.

\subsection{Existence of blenders when at least one of $\theta$ and $\dfrac{\omega}{2\pi}$ is irrational}\label{sec:existsf}
Recall that for a saddle-focus heterodimensional cycle, we denote the complex center-stable multipliers by $\lambda e^{\pm i\omega}$, where $\lambda>0$ and $0<\omega<\pi$, and the real center-unstable multiplier by $\gamma$. Let $A>0,B>0,b,u^-,\eta_{1,2}$ be the coefficients in \eqref{akmbkm}. We introduce a condition on the modulus $\omega$. 

\noindent\setword{\textbf{A1}}{word:A1}: there exist infinitely many $k>0$ such that
\begin{equation}\label{eq:bad}
bu^-\sin(k\omega +\eta_2)>0,\quad \sin(k\omega +\eta_1)\neq 0,\quad 
\alpha_k:=\dfrac{bu^-A\sin(k\omega+\eta_1)}{\sin(k\omega+\eta_2)}\neq\pm 1.
\end{equation}

Observe that condition \ref{word:A1} is automatically satisfied if $\omega/2\pi$ is irrational, or if $\omega/2\pi = p/q$ with $q\geqslant 9$. To see the latter, note first that the last inequality in \eqref{eq:bad} can be rewritten as
$$|bu^- \sin(k\omega +\eta_1)|\neq|\sin(k\omega +\eta_1)\cos(\eta_2-\eta_1)
+\sin(k\omega +\eta_1)\sin(\eta_2-\eta_1)|,$$
or
$$|(bu^--\cos(\eta_2-\eta_1))\tan(k\omega +\eta_1)|\neq|\sin(\eta_2-\eta_1)|.$$
It follows that to achieve simultaneously this inequality and the second one in \eqref{eq:bad} it is enough to have at least 4 different values of $k\omega\bmod 2\pi$  in $[0,1/2)$ or in $[1/2,1)$. It is easy to check that, when $\omega/2\pi=p/q$, there are at least $\floor{(q-1)/2}$ different values of $k\omega\bmod 2\pi$ both in $[0,1/2)$ and in $[1/2,0)$, where $\floor{a}$ means the largest integer which is smaller than or equal to $a$. Choosing the  values from the suitable interval also gives $bu^-\sin(k\omega +\eta_2)>0$. Thus, when $\floor{(q-1)/2}\geqslant 4$ or $q\geqslant 9$, there  exists some $r\in \{0,\dots,q-1\}$ such that all $k=r+iq$ with $i=0,1,2,\dots$ satisfy  the inequalities in \eqref{eq:bad}.

\begin{prop}\label{prop:blendersf}
Let $\delta$ be the size of $\Pi$ defined in \eqref{eq:complex:domain}, and $\mathcal{C}^{ss}$ and $\mathcal{C}^{uu}$ be given by Lemma \ref{lem:conefields} with a sufficiently small $K$. Let a non-degenerate saddle-focus heterodimensional cycle satisfy condition \ref{word:A1} and $U$ be any neighborhood of the cycle. 
\begin{enumerate}[nosep]
\item If $\theta=\ln \lambda/\ln |\gamma|$ is irrational and $\omega/2\pi=p/q$, then, for every sufficiently large $k$ satisfying condition \ref{word:A1}, there exists a cs-blender in $U$  when $|\alpha_k|<1$ (defined in \eqref{eq:bad}) or a cu-blender in $U$ when $|\alpha_k|>1$.

\item If $\omega/2\pi$ is irrational, and $\theta,\omega/2\pi$ and 1 are not rationally independent, then there exists a cu-blender in $U$.

\item If $\theta,\omega/2\pi$ and 1 are rationally independent, then there exist simultaneously  a cs-blender and a cu-blender in $U$.
\end{enumerate}
In each of the above cases, the cs-blender has index $d_1$ and an activating pair $(\Pi',\mathcal{C}^{ss})$, while  the cu-blender has index $(d_2=d_1+1)$ and an activating pair $(\Pi',\mathcal{C}^{uu})$, where the cube $\Pi'$ is defined as
\begin{equation}\label{eq:Pi'sf}
\Pi'=[-c\delta,c\delta]\times[-\delta,\delta]^{d_1}\times [-\delta,\delta]^{d-d_1-1}\subset \Pi,
\end{equation}
for some $c\in(0,1)$.
\end{prop}

\begin{proof} We prove the proposition by finding in each case a  sequence of pairs $\{(k_n,m_n)\}$  such that the sequence of return maps $\{T_{k_n,m_n}\}$ satisfies the conditions of Proposition \ref{prop:blender}. Let $A>0,B>0,b,u^-,\eta_{1,2}$ be the coefficients in \eqref{eq:complex:T_k,m_cross_0mu}.

\noindent\textbf{Case 1: $\theta$ is irrational and $\omega/2\pi$ is rational.}
 Take $\omega/2\pi=p/q$. Since condition \ref{word:A1} is satisfied, there exists $k=r\in \{0,\dots,q-1\}$ satisfying \eqref{eq:bad}. We particularly have
\begin{equation}\label{eq:sf:r}
\sin \left( \dfrac{2\pi rp}{q} +\eta_2\right)Bbu^-=: R Bbu^->0.
\end{equation}
Since $\theta$ is irrational, the set 
\begin{equation*}\label{eq:sf:set1}
\{m-(kq+r)\theta \}_{k\in \mathbb{N},m\in 2\mathbb N}
\end{equation*}
is dense in $\mathbb R$, and, hence, $\{\gamma^{m-(kq+r)\theta }\}_{k\in \mathbb{N},m\in 2\mathbb N}$ is dense in $[0,\infty]$. This along with \eqref{eq:sf:r} implies that one can find a sequence $\{(k'_nq+r,m'_n)\}$  such that $\{\gamma^{m'_n-(k'_nq+r)\theta }\}$ is dense in a small neighborhood of $bu^-/BR$. Thus, the sequence
\begin{equation}\label{eq:sf:sequence}
\{|\gamma|^{m'_n-(k'_nq+r)\theta }B\sin((k'_nq+r)\omega+\eta_2)-bu^-\}
\end{equation}
is dense in a small neighborhood $\Delta$ of zero.

Recall $\theta=-\ln\lambda/\ln|\gamma|$ (note that $\lambda>0$). When the integer $m$ is even, we can write
\begin{equation}\label{eq:lambdagamma}
\lambda^k\gamma^m=|\gamma|^{m-\theta k}.
\end{equation}
Now denote $k_n=k'_nq+r$ and $m_n=m'_n$. The preceding paragraph shows that $\sin(k_n\omega+\eta_2)$ is bounded away from zero. Also, since $\lambda>0$ and $m_n$ are even, relation \eqref{eq:lambdagamma} implies that $\lambda^{k_n}\gamma^{m_n}=\gamma^{m_n-k_n\theta}$. It follows that by taking $\Delta$ sufficiently small, the condition \eqref{eq:crosscon} is satisfies with $(k,m)=(k_n,m_n)$. Hence, the first return maps $T_n:=T_{k_n,m_n}$  take the form \eqref{eq:complex:T_k,m_cross_0mu}. By comparing \eqref{akmbkm} with \eqref{eq:sf:sequence} we find that the sequence \eqref{eq:sf:sequence} is exactly the constant term $B_{k_n,m_n}$ in \eqref{eq:complex:T_k,m_cross_0mu}. Thus, setting $X=X_1$ and $Z^{new}=(X_2,Z)$, we find that  the sequence $\{T_n\}$ satisfies condition 1 of Proposition \ref{prop:blender}.

Next, let us check condition 2 of Proposition \ref{prop:blender}. To find the values of $|A_{k_n,m_n}|$, we use $\eqref{akmbkm}$ to write
\begin{equation}\label{eq:sf:Akm}
A_{k_n,m_n}=\dfrac{(B_{k_n,m_n}+bu^-)A\sin(k_n\omega+\eta_1)}{\sin(k_n\omega+\eta_2)}.
\end{equation}
On the one hand, by taking $\Delta$ sufficiently small the contribution of $B_{k_n,m_n}$ to $A_{k_n,m_n}$ can be ignored. On the other hand, by the choice of $k_n$ one has $\sin(k_n\omega+\eta_1)=\sin(r\omega+\eta_1)$ and $\sin(k_n\omega+\eta_2)=\sin(r\omega+\eta_2)$. Thus, we have
\begin{equation}\label{eq:sf:Akm2}
A_{k_n,m_n}\sim \dfrac{bu^-A\sin(r\omega+\eta_1)}{\sin(r\omega+\eta_2)}=:\alpha_r,
\end{equation}
Since assumption  \eqref{eq:bad} ensures $|\alpha_r|\neq 1$, this gives  condition 2.1 of Proposition \ref{prop:blender}  if $|\alpha_r|<1$, or condition 2.2 if $|\alpha_r|>1$. Therefore, using Proposition \ref{prop:blender} we obtain in the former case a cs-blender with an activating pair $(\Pi',\mathcal{C}_1)$ for some cone field $\mathcal{C}_1$, and in the latter case a cu-blender with an activating pair $(\Pi',\mathcal{C}_2)$ for some cone field $\mathcal{C}_2$. 

Since these two cone fields are independent of $\delta$ and $n$, one can  take $\delta$ small and $n$ large so that  $\mathcal{C}^{ss}$ and $\mathcal{C}^{uu}$ in Lemma \ref{lem:conefields} are contained in $\mathcal{C}_1$ and $\mathcal{C}_2$, respective. As a result, $(\Pi',\mathcal{C}^{ss})$ and $(\Pi',\mathcal{C}^{uu})$ are also activating pairs.


\noindent\textbf{Case 2: $\omega/2\pi$ is irrational, and $\theta,\omega/2\pi$ and 1 are rationally dependent.} We can write
\begin{equation*}\label{eq:rtdp} 
\theta =\dfrac{p}{q} \dfrac{\omega}{2\pi} + \dfrac{p'}{q'}.
\end{equation*}
Consider the function 
$$h(s)=|\gamma|^{ -\frac{p}{q}s}B\sin (2\pi s +\eta_2)-bu^-,$$
which is defined on $\mathbb{R}$. Let  $s=s_0$ be a zero of $h$ , i.e., $h(s_0)=0$. It is obvious that $s_0$ can be taken with arbitrarily large absolute value.
 Since the set 
\begin{equation*}\label{eq:sf:set2.2}
\left\{qq'k\dfrac{\omega}{2\pi}-qq'l\right\}_{k,l\in 2\mathbb{N}}
\end{equation*}
is dense in $\mathbb R$, one can find a sequence $\{qq'k'_n\omega/2\pi-qq'l_n\}$ dense in a small neighborhood of $s_0$ such that, by continuity of $h(s)$, the sequence
\begin{equation}\label{eq:sf:2}
\{|\gamma|^{m_0 -\frac{p}{q}(qq'k'_n\frac{\omega}{2\pi}-qq'l_n)}B\sin (qq'k'_n \omega +\eta_2)-bu^-\}
\end{equation}
is dense in a small neighborhood $\Delta$  of zero. 

By \eqref{eq:lambdagamma} and \eqref{eq:rtdp}, for even $m$ we have 
$$\lambda^k\gamma^m=|\gamma|^{m-k\frac{p}{q}\frac{\omega}{2\pi}-k\frac{p'}{q'}}.$$
Taking $k_n=qq'k'_n$ and $m_n=p'qk'_n+pq'l_n$, one readily finds that $B_{k_n,m_n}$ in \eqref{eq:complex:T_k,m_cross_0mu} equals \eqref{eq:sf:2}. Arguing as in case (1), we obtain condition 1 of Proposition \ref{prop:blender} with $\{T_n:=T_{k_n,m_n}\}$.

It is immediate from \eqref{eq:sf:2} that  $\sin(k_n\omega+\eta_2)$ can be arbitrarily close to zero by taking  $|s_0|$ sufficiently large. Since $B_{k_n,m_n}$ is negligible by taking $\Delta$ small, it  follows from \eqref{eq:eta} and \eqref{eq:sf:Akm} that $|A_{k_n,m_n}|>1$ for all $n$, and, particularly, they are bounded away from $1$ and infinity. So, applying Proposition \ref{prop:blender} with condition 2.2 gives a cu-blender. The activating pair is obtained as in case 1.

\noindent\textbf{Case 3: $\theta,\omega/2\pi$ and 1 are rationally independent.} The result is an easy consequence of the denseness of the set 
\begin{equation*}\label{eq:sf:set3}
\left\{\left(k\theta -m,k\dfrac{\omega}{2\pi} - l\right)\right\}_{k,\in\mathbb N,m\in2\mathbb N,l\in\mathbb Z}
\end{equation*}
in $\mathbb R^2$. Indeed, one just needs to choose a sequence $\{(k_n,m_n)\}$ satisfying  that $\sin(k_n\omega+\eta_2)$ is dense in a small interval $\Delta'$ and $\lambda^{k_n}\gamma^{m_n}=|\gamma|^{m_n-k_n\theta} \to t$ for some constant $t$ such that
\begin{equation*}\label{eq:sf:3}
\{B_{k_n,m_n}=|\gamma|^{m_n-k_n\theta }B\sin(k_n\omega+\eta_2)-bu^-\}
\end{equation*}
is dense in a small neighborhood $\Delta$ of zero. The latter fact  gives condition 1 of Proposition \ref{prop:blender}. The former one implies that, with sufficiently small $\Delta$, the values $|A_{k_n,m_n}|$ in \eqref{eq:sf:Akm} can lie between zero and 1 by taking $\Delta'$ close to zero, which leads to a cs-blender, or lie between 1 and infinity by taking $\Delta'$ such that $\sin(k_n\omega+\eta_2)$ is close to zero, which leads to a cu-blender. The activating pair is obtained as in case 1.
\end{proof}


\subsection{Homoclinic relations in the saddle-focus case}\label{sec:homosf}
We first make some observations on the   blenders obtained in Section \ref{sec:existsf}. The cs-blender in Proposition~\ref{prop:blendersf}, denoted by $\Lambda_1$, is given by Proposition~\ref{prop:blender}, and hence corresponds to a finite set $\mathcal{P}_1$ of pairs $(k,m)$ such that $0<C_1<|A_{k,m}|<C_2<1$ for $(k,m)\in \mathcal{P}_1$, where $A_{k,m}$ is as in \eqref{eq:complex:T_k,m_cross_0mu} and $C_{1,2}$ are some constants. Similarly, a cu-blender in Proposition~\ref{prop:blendersf}, denoted by $\Lambda_2$, corresponds to a finite set $\mathcal{P}_2$ of pairs $(k,m)$ such that $1<C_3<|A_{k,m}|<C_4<\infty$ for $(k,m)\in \mathcal{P}_2$, where $C_{3,4}$ are some constants. The sets $\Lambda_i$ $(i=1,2)$ consist of points whose orbits under $T_{k,m}$ with $(k,m)\in \mathcal{P}_i$ never leave the cube $\Pi$ defined in \eqref{eq:complex:domain} and hence are hyperbolic basic sets equipped with the cone fields in Lemma \ref{lem:conefields}.

For each $i=1,2$ and a point $M_i\in \Lambda_i\cap\Pi'$, define its local invariant manifolds $W^s_{loc}(M_i)$ and $W^u_{loc}(M_i)$ as the connected components through $M_i$ of $W^s_{loc}(\Lambda_i)\cap\Pi$ and, respectively, of $W^u_{loc}(\Lambda_i)\cap\Pi$. One readily finds from Lemma \ref{lem:conefields} that those manifolds are given by   smooth functions of the following forms:
\begin{itemize}[nosep]
\item $W^s_{loc}(M_1):Y=\psi^s_{M_1}(X,Z)$ defined on $[-\delta,\delta]\times[-\delta,\delta]^{d-d_1-1}$;
\item $W^u_{loc}(M_1):(X,Z)=\psi^u_{M_1}(Y)$ defined on $[-\delta,\delta]^{d_1}$, and, in particular, $W^u_{loc}(M_1)$ crosses $\Pi'$ properly with respect to $\mathcal{C}^{uu}$ defined in \eqref{eq:uucone};
\item $W^s_{loc}(M_2):(X,Y)=\psi^s_{M_2}(Z)$ defined on $[-\delta,\delta]^{d-d_1-1}$, and, in particular, $W^s_{loc}(M_2)$ crosses $\Pi'$ properly with respect to $\mathcal{C}^{ss}$ defined in \eqref{eq:sscone};
\item $W^u_{loc}(M_2):Z=\psi^u_{M_2}(X,Y)$ defined on $[-\delta,\delta]\times[-\delta,\delta]^{d_1}$.
\end{itemize}

An immediate observation is that when a cs-blender and a cu-blender found in Proposition \ref{prop:blendersf} coexist (i.e.,   case 3), they activate each other in the sense of Definition \ref{defi:act}.


%

For the other cases in Proposition~\ref{prop:blendersf}, we need the fact (see equation (3.23) of \citep{LT21}) that
\begin{equation}\label{eq:Omf}
W^s_{loc}(O_1)\cap\Pi=\{Y=0\}\quad\mbox{and}\quad
F_{21}(W^u_{loc}(O_2))\cap\Pi=\{Z=0\},
\end{equation}
and the following 
\begin{lem}[{\citep[Lemma 6.5]{LT21}}] \label{lem:LTlemmasf}
If $\omega/2\pi$ is irrational, then there exist infinitely many $k$ such that $F^{-k}_1\circ F^{-1}_{12}(W^s_{loc}(O_2))$ is a disc of the form
\begin{equation*}\label{eq:Omf2}
(X,Y)=s(Z),
\end{equation*}
which crosses $\Pi'$ defined by \eqref{eq:Pi'sf} properly with respect to $\mathcal{C}^{ss}$ defined by \eqref{eq:sscone}. 
\end{lem}

Comparing the equations in \eqref{eq:Omf}  with those for $W^u_{loc}(M_1)$ and $W^s_{loc}(M_2)$, we obtain 
\begin{equation}\label{eq:autoint1}
W^u(O_2)\pitchfork W^s(\Lambda_2)\neq \emptyset
\quad\mbox{and}\quad
W^s(O_1)\pitchfork W^u(\Lambda_1)\neq \emptyset.
\end{equation}
Here the second intersection will be used later in the proof of the perturbative results. By the first intersection, the formulas  for $W^s_{loc}(M_2)$ and $W^u_{loc}(M_2)$, and Lemma~\ref{lem:LTlemmasf}, one  finds that $O_2$ is homoclinically related to the cu-blender found in Proposition \ref{prop:blendersf} when $\omega/2\pi$ is irrational (i.e. cases 1 and 2). In this case, it follows from the lambda-lemma that $O_2$ also activates the cs-blender if the cu-blender does. In summary, we obtain

\begin{lem}\label{lem:homo_sf}
The periodic point $O_{2}$ and the blenders found in Proposition \ref{prop:blendersf} have the following homoclinic relations:
\begin{itemize}[nosep]
\item in case 2, the cu-blender is homoclinically related to $O_2$, and
\item in case 3, the cu-blender is homoclinically related to $O_2$; the two blenders activate each other, and in particular $O_2$ activates the cs-blender.
\end{itemize}
\end{lem}

\subsection{Proofs of Theorem \ref{thm:sf} and the first case in Theorem \ref{thm:rational2}}\label{sec:proofsf}

Theorem \ref{thm:sf} follows immediately from Proposition \ref{prop:blendersf} and Lemma \ref{lem:homo_sf}.

\begin{proof}[Proof of the first case in Theorem~\ref{thm:rational2}]
Since being non-degenerate is an open and dense property in $\mathcal{H}_f$, we take the neighborhood $\mathcal{U}$ of $g$ sufficiently small such that $\Gamma_f$ is non-degenerate for all $f\in\mathcal{U}$.
Take any $f\in\mathcal{U}$ with irrational $\theta(f)$ and rational $\omega(f)/2\pi=p/q$. According to the proof of Proposition~\ref{prop:blender},  to obtain both a cs-blender and a cu-blender for $f$, it suffices to find $r_1,r_2\in\{0,\dots,q-1\}$ such that $0<|\alpha_{r_1}|<1$ and $1<|\alpha_{r_2}|<\infty$ (see  \eqref{eq:sf:Akm2}). Obviously, there exists $N_f$ depending on $A,b,u^-,\eta_{1,2}$ such that if $q\geq N_f$  then one can always find such $r_1$ and $r_2$.  Since the cycle is non-degenerate for all $f\in\mathcal{U}$,  the involved coefficients are well defined and depend continuously on $f$. Thus, coexisting blenders appear whenever $q\geq N(\mathcal{U})=\sup_{f\in\mathcal{U}}N_f$. The homoclinic relation follows from the formulas for the invariant manifolds of points of the blenders listed in Section~\ref{sec:homosf}.
\end{proof}

\section{Blenders near double-focus heterodimensional cycles}\label{sec:blenderdf}
In this section, we prove Theorem \ref{thm:df} and the second case in Theorem \ref{thm:rational2}, in the same fashion as for saddle-focus cycles. We first state the blender results and homoclinic relation results separately, and then prove the Theorems.

\subsection{Existence of blenders when at least one of $\theta,\dfrac{\omega_1}{2\pi},\dfrac{\omega_2}{2\pi}$ is irrational}
Recall that for a double-focus heterodimensional cycle, we denote the  center-stable multipliers by $\lambda e^{\pm i\omega_1}$ and the center-unstable multipliers by $\gamma e^{\pm i\omega_2}$, where $\lambda>0,\gamma>0$ and $0<\omega_{1,2}<\pi$.  Let $A^*>0,B^*>0,b,u^-,\eta_{1,2,3}$ be the coefficients in \eqref{eq:df:AkmBkm} and \eqref{eq:short1}. We first introduce a condition on the moduli $\omega_1$ and $\omega_2$.

\noindent\setword{\textbf{A2}}{word:A2}:  
there exist infinitely many  pairs $(k,m)$  satisfying \eqref{eq:quant} and the inequalities
\begin{equation}\label{eq:bad2}
\begin{aligned}
&(\cos(m\omega_2+\eta_3)u^-_1+\sin(m\omega_2+\eta_3)u^-_2)\sin(k\omega_1+\eta_2)>0,
\qquad \sin(k\omega_1 +\eta_1)\neq 0,\\[5pt]
&\alpha_{k,m}:=\dfrac{A^*(\cos(m\omega_2+\eta_3)u^-_1+\sin(m\omega_2+\eta_3)u^-_2)}{(\cos(m\omega_2+\eta_3) + b\sin(m\omega_2+\eta_3))^2}
\cdot 
\dfrac{\sin(k\omega_1+\eta_1)}{\sin(k\omega_1+\eta_2)}\neq\pm 1.
\end{aligned}
\end{equation}

Obviously, condition \ref{word:A2} is satisfied when both $\omega_1/2\pi$ and $\omega_2/2\pi$ are irrational. This condition also holds  if $\omega_1/2\pi$ is irrational and $\omega_2/2\pi=p/q$ with $q\geqslant \q2$, and if $\omega_1/2\pi=p_1/q_1$ with $q_1\geqslant 9$ and $\omega_2/2\pi=p_2/q_2$ with $q_2\geqslant \q2$. Indeed, when $q_2\geq \q2$, there are at least 6 different values of $m\omega_2\bmod 2\pi$ in $[0,2\pi)$, so at least one value will satisfy \eqref{eq:quant} and $\cos(m\omega_2+\eta_3)u^-_1+\sin(m\omega_2+\eta_3)u^-_2\neq 0$. Taking this value fixes the factors involving $m$ in   \eqref{eq:bad2}. Then,  arguing as in the saddle-focus case for condition \ref{word:A1}, we find that, when $\omega_1$ is irrational  or $q_1\geq 9$, there is at least one value of $k\omega_1 \bmod 2\pi$ validating all inequalities in \eqref{eq:bad2}. 

For the sake of clarity in proving our results, we classify the various possibilities of the arithmetic relations among $\theta,\omega_1/2\pi,\omega_2/2\pi$ into 7 cases.

\begin{prop}\label{prop:blenderdf}
Let $\delta$ be the size of $\Pi$ defined in \eqref{eq:df:domain}, and $\mathcal{C}^{ss}$ and $\mathcal{C}^{uu}$ be given by Lemma \ref{lem:conefields_df} with a sufficiently small $K$. Let a non-degenerate double-focus heterodimensional cycle satisfy condition~\ref{word:A2} and $U$ be any neighborhood of the cycle. 
\begin{enumerate}[nosep]
\item[1.] If $\theta$ is irrational, and $\omega_1/2\pi,\omega_2/2\pi$ are rational, then, for every sufficiently large pair $(k,m)$ satisfying condition \ref{word:A2}, there exists  a cs-blender in $U$ if $|\alpha_{k,m}|<1$  (defined in \eqref{eq:bad2}), or a cu-blender in $U$ if $|\alpha_{k,m}|>1$.
\item[2.1.] If $\omega_1/2\pi$ is irrational, $\omega_2/2\pi$ is rational, and $\theta,\omega_1/2\pi,1$ are rationally independent, then there exist simultaneously a cs-blender and a cu-blender in $U$.
\item[2.2.] If  $\omega_1/2\pi$ is irrational, $\omega_2/2\pi$ is rational, and $\theta,\omega_1/2\pi,1$ are rationally dependent, then there exists a cu-blender in $U$.
\item[3.1.] If $\omega_1/2\pi,\omega_2/2\pi$ are irrational, and $\theta,\omega_1/2\pi,\theta\omega_2/2\pi,1$ are rationally independent, then there exist simultaneously a cs-blender and a cu-blender in $U$.
\item[3.2.] If $\omega_1/2\pi,\omega_2/2\pi$ are irrational, $\theta,\omega_1/2\pi,\theta\omega_2/2\pi,1$ are rationally dependent, but $\theta,\omega_1/2\pi,1$ are rationally independent, then there exist simultaneously a cs-blender and a cu-blender in $U$.
\item[3.3.1.] If $\omega_1/2\pi,\omega_2/2\pi$ are irrational, $\theta,\omega_1/2\pi,1$ are rationally dependent, and $\omega_1/2\pi,\theta\omega_2/2\pi,1$ are rationally dependent, then there exists a cu-blender in $U$.
\item[3.3.2.] If $\omega_1/2\pi,\omega_2/2\pi$ are irrational, $\theta,\omega_1/2\pi,1$ are rationally dependent, and $\omega_1/2\pi,\theta\omega_2/2\pi,1$ are rationally independent, then there exists a cu-blender in $U$.
\end{enumerate}
In each of the above cases, the cs-blender has index $d_1$ and an activating pair $(\Pi',\mathcal{C}^{ss})$, while  the cu-blender has index $(d_2=d_1+1)$ and an activating pair $(\Pi',\mathcal{C}^{uu})$, where the cube $\Pi'$ is defined as
\begin{equation}\label{eq:Pi'df}
\Pi'=[-c\delta,c\delta]\times[-\delta,\delta]^{d-d_2}\times [-\delta,\delta]^{d_2-1}\subset \Pi,
\end{equation}
for some $c\in(0,1)$.
\end{prop}

We postpone the proof of this proposition to Section \ref{sec:proofdf}, after we prove the theorems.

\subsection{Homoclinic relations in the double-focus case}\label{sec:homodf}
Let us find out the homoclinic relations among the found blenders in Proposition~\ref{prop:blenderdf} and the periodic points $O_1$ and $O_2$. Denote by $\Lambda_1$ and $\Lambda_2$ the cs- and cu-blenders.  Since the proof of Proposition~\ref{prop:blenderdf} (see Section~\ref{sec:proofprop}) is also based on the use of Proposition \ref{prop:blender}, the structures of $\Lambda_i$ $(i=1,2)$  are completely the same as those in the saddle-focus case, i.e., they consist of points  whose orbits under the corresponding sets of first return maps never leave the cube $\Pi$ defined in \eqref{eq:df:domain}, and their activating domain is $\Pi'$ defined in \eqref{eq:Pi'df}.
 (See the discussion at the beginning of Section~\ref{sec:homosf}). 
One concludes from Lemma~\ref{lem:conefields_df} that for any points $M_1\in \Lambda_1\cap\Pi'$ and $M_2\in \Lambda_2\cap\Pi'$, their local invariant manifolds  are given by   smooth functions of the following forms:
\begin{itemize}[nosep]
\item $W^s_{loc}(M_1):W=\psi^s_{M_1}(U,V)$ defined on $[-\delta,\delta]\times[-\delta,\delta]^{d-d_2}$ with $c$ in \eqref{eq:Pi'df};
\item $W^u_{loc}(M_1):(U,V)=\psi^u_{M_1}(W)$ defined on $[-\delta,\delta]^{d_2-1}$, and, in particular, $W^u_{loc}(M_1)$ crosses $\Pi$ properly with respect to $\mathcal{C}^{uu}$ defined in \eqref{eq:uucone_df};
\item $W^s_{loc}(M_2):(U,W)=\psi^s_{M_2}(V)$ defined on $[-\delta,\delta]^{d-d_2}$, and, in particular, $W^s_{loc}(M_2)$ crosses $\Pi$ properly with respect to $\mathcal{C}^{ss}$ defined in \eqref{eq:sscone_df};
\item $W^u_{loc}(M_2):V=\psi^u_{M_2}(U,W)$ defined on $[-\delta,\delta]\times [-\delta,\delta]^{d_2-1}$.
\end{itemize}

Obviously, whenever two blenders in Proposition~\ref{prop:blenderdf} coexist (i.e., cases 2.1, 3.1 and 3.2), they activate each other in the sense of Definition~\ref{defi:act}.

For other cases in Proposition~\ref{prop:blenderdf}, we need the fact, given at the end of \citep[Section 6.4]{LT21}, that
$$
W^u_{loc}(O_2)=\{V=0\}
\quad\mbox{and}\quad
F^{-1}_{21}(W^s_{loc}(O_1))\cap \Pi'=\{W=\psi(U,V)\},
$$
where $\psi=O(U^2+\delta V^2)$ is some smooth function defined on $[-c\delta,c\delta]\times[-\delta,\delta]^{d-d_2}$. Comparing these equations with those for $W^s_{loc}(M_2)$ and $W^u_{loc}(M_1)$, we obtain that 
\begin{equation}\label{eq:autoint2}
W^u(O_2)\pitchfork W^s(\Lambda_2)\neq \emptyset
\quad\mbox{and}\quad
W^s(O_1)\pitchfork W^u(\Lambda_1)\neq \emptyset.
\end{equation}
For other heteroclinic intersections, we use the following

\begin{lem}[{\citep[Lemma 6.8]{LT21}}] \label{lem:LTlemma14}
If $\omega_2/2\pi$ is irrational, then there exist infinitely many $m$ such that $F^{m}_2\circ F_{12}(W^u_{loc}(O_1))\cap \Pi'$ is a disc of the form 
$(U,V)=s(W),$
crossing $\Pi'$ properly with respect to $\mathcal{C}^{uu}$ defined by \eqref{eq:uucone_df}.
\end{lem}

\begin{lem}[{\citep[Lemma 6.9]{LT21}}] \label{lem:LTlemma15}
If $\omega_1/2\pi$ is irrational, then there exist infinitely many $k$ such that $F^{-1}_{21}\circ F^{-k}_1\circ F^{-1}_{12}(W^s_{loc}(O_2))$ is a disc of the form 
$(U,W)=s_2(V),$
crossing $\Pi'$ defined by \eqref{eq:Pi'df} properly with respect to $\mathcal{C}^{ss}$ defined by \eqref{eq:sscone_df}.
\end{lem}

Thus, if $\omega_2/2\pi$ is irrational, then $W^u(O_1)\pitchfork W^s(\Lambda_1)\neq \emptyset$, which along with \eqref{eq:autoint2} implies that  $O_1$ is homoclinically related to the cs-blender $\Lambda_1$;  if $\omega_1/2\pi$ is irrational, then $W^s(O_2)\pitchfork W^u(\Lambda_2)\neq \emptyset$, which along with \eqref{eq:autoint2} implies that  hence $O_2$ is homoclinically related to the cu-blender $\Lambda_2$. In summary, we obtain

\begin{lem}\label{lem:homo_df}
The periodic points $O_{1,2}$ and the blenders found in Proposition \ref{prop:blenderdf} have the following homoclinic relations:
\begin{itemize}[nosep]
\item in case 2.1, the cu-blender is homoclinically related to $O_2$; the two blenders activate each other, and, in particular, $O_2$ activates the cs-blender,
\item in case 2.2, the cu-blender is homoclinically related to $O_2$, 
\item in cases 3.1 and 3.2, the cs-blender is homoclinically related to $O_1$ and the cu-blender is  homoclinically related to $O_2$; the two blenders activate each other, and, in particular, $O_1$ activates the cu-blender and  $O_2$ activates the cs-blender, and
\item in cases 3.3.1 and 3.3.2, the cu-blender is homoclinically related to $O_2$ and activated by $O_1$.
\end{itemize}
\end{lem}

\subsection{Proofs of Theorem \ref{thm:df} and the second case in Theorem \ref{thm:rational2}}\label{sec:proofdf}

\begin{proof}[Proof of Theorem \ref{thm:df}.]
Cases 2.1 and 2.2 in Proposition \ref{prop:blenderdf} together with Lemma \ref{lem:homo_df} cover case (1.i), and, by considering $f^{-1}$, also case (1.ii). Since the rational independence of $\theta,\omega_1/2\pi,1$ implies the irrationality of $\theta$, case (2.i) follows from case 2.1 in Proposition \ref{prop:blenderdf} and Lemma \ref{lem:homo_df}. Using $f^{-1}$, we also obtain case (2.ii). Case (3) follows immediately from case 1 in Proposition \ref{prop:blenderdf}.

What remains is case (1.iii), whose hypothesis is covered by those of cases 3.1, 3.2, 3.3.1 and 3.3.2 in Proposition \ref{prop:blenderdf}. We are done if either of the hypotheses in cases 3.1 and 3.2  holds. Now suppose we are under the hypothesis of case 3.3.1. With Lemma \ref{lem:homo_df}, this gives us a cu-blender homoclinically related to $O_2$ and activated by $O_1$. Considering $f^{-1}$ leads  to  one of cases 3.1, 3.2, 3.3.1 and 3.3.2. Since the roles of $O_1$ and $O_2$ interchange for $f^{-1}$, it follows from  Proposition \ref{prop:blenderdf} and Lemma \ref{lem:homo_df} that we always have a cu-blender of $f^{-1}$, homoclinically related to $O_1$ and activated by $O_2$. This blender is a cs-blender of $f$. Using the lambda-lemma, we see that this cs-blender and the previous cu-blender activate each other, giving the desired result for case (1.iii). The proof with the hypothesis of case 3.3.1 is the same.
\end{proof}

\begin{proof}[Proof of the second case in Theorem \ref{thm:rational2}]
The proof is essentially the same as in the saddle-focus case. The only difference is  that the quantity $\alpha_k$ there is   replaced by $\alpha_{k,m}$ defined in \eqref{eq:bad2}, which depends on the coefficients $A^*,b,u^-_1,u^-_2,\eta_{1,2,3}$ in \eqref{eq:df:AkmBkm}.
Take the neighborhood $\mathcal{U}$ of $g$ sufficiently small such that $\Gamma_f$ is non-degenerate for all $f\in\mathcal{U}$.
 Take any $f\in\mathcal{U}$ with rational $\omega_1(f)/2\pi$ and $\omega_2(f)/2\pi$, and denote their denominators by $q_1$ and $q_2$, respectively. Since the coefficients depend on $f$ continuously, there exists a lower bound $N_1(\mathcal{U})$ such that if $q_1\geq N_1(\mathcal{U})$ and $q_2\geq \q2$, then one can find $r_1,r_1'\in \{0,\dots,q_1-1\}$ and $r_2\in \{0,\dots,q_2-1\}$ satisfying that $|\alpha_{r_1,r_2}|<1$ and $|\alpha_{r'_1,r_2}|>1$. Note that here taking $q_2\geq \q2$ is only for the validation of condition \ref{word:A2}, and after fixing $r_2$ one only needs to change $r_1$ adjust the values of $|\alpha_{r_1,r_2}|$. Applying the same arguments to $f^{-1}$ gives a lower bound $N_2(\mathcal{U})$ for $q_2$. Setting $N(\mathcal{U})=\max \{N_1(\mathcal{U}),N_2(\mathcal{U})\}$ proves the theorem.
\end{proof}

\subsection{Proof of Proposition \ref{prop:blenderdf}}\label{sec:proofprop}
We first prove an auxiliary result (Lemma \ref{lem:reduction} below). 
Define two functions
\begin{align}
&g:\mathbb{R}\times S^1\times S^1\to \mathbb R,\quad g(t,s,w)=\gamma^tA(2\pi  w)\sin(2\pi s+\eta_1), \label{eq:g}\\
&h:\mathbb{R}\times S^1\times S^1\to \mathbb R,\quad h(t,s,w)=\gamma^tB(2\pi  w)\sin (2\pi  s+\eta_2) - \xi(2\pi w), \label{eq:h}
\end{align}
where $S^1=[0,1]$
with the endpoints identified, and $A,B,\xi$ are given in \eqref{eq:short1}. We consider two classes of lines and one class of planes in $\mathbb{R}^3$: 
\begin{flalign}\label{eq:line}
\text{(proper horizontal line)} &&\{s=s^*,w=w^*\},&&
\end{flalign}
where $s^*$ and $w^*$ are constants  satisfying 
\begin{equation}\label{eq:con_line}
0<\xi(2\pi w^*)B(2\pi w^*)\sin(2\pi  s^*+\eta_2)<\infty,\quad
\dfrac{\xi(2\pi w^*)A(2\pi w^*)\sin(2\pi  s^*+\eta_1)}{\sin(2\pi  s^*+\eta_2)}\not\in\{0,\pm 1,\infty\};
\end{equation}
\begin{flalign}\label{eq:segment}
\text{(proper tilted line)} &&\{t=a_1s+b,w=a_2s\}\quad\mbox{or}\quad \{t=a_1s+b,w=w^*\},&&
\end{flalign}
where $a_1^2+b^2\neq 0,a_2\neq 0$, and $w^*$ satisfies $A(2\pi w^*)B(2\pi w^*)\xi(2\pi w^*)\not\in \{0,\infty\}$;
\begin{flalign}\label{eq:plane}
\text{(proper plane)} &&\{w=at+bs\}\quad\mbox{or}\quad \{w=w^*\},&&
\end{flalign}
where $a^2+b^2\neq 0$, and $w^*$ satisfies $A(2\pi w^*)B(2\pi w^*)\xi(2\pi w^*)\not\in \{0,\infty\}$.

  
\begin{lem}\label{lem:reduction}
For a sequence $E=\{(t_n,s_n,w_n)\}$ of points in $\mathbb{R}^3$,
\begin{enumerate}[nosep]
\item if the closure $cl(E)$ contains a proper horizontal line, then there exists a subsequence for which $\{h(t_n,s_n,w_n)\}$ is dense near zero and $\{|g(t_n,s_n,w_n)|\}$ belongs to  some closed interval either in $(0,1)$ or in $(1,\infty)$, depending on $s^*$ and $w^*$ in \eqref{eq:line} and coefficients in $g$ and $h$;
\item if the closure $cl(E)$ contains a proper tilted line, then there exists a subsequence for which $\{h(t_n,s_n,w_n)\}$ is dense near zero and $\{|g(t_n,s_n,w_n)|\}$ belongs to  some closed interval in $(1,\infty)$;
\item if the closure $cl(E)$ contains a proper plane, then there exist two subsequences such that $\{h(t_n,s_n,w_n)\}$ is dense near zero for both of them and $\{|g(t_n,s_n,w_n)|\}$ belongs to some closed interval in $(0,1)$ for one of them and in $(1,\infty)$ for the other.
\end{enumerate}
\end{lem}
\begin{proof}
To simplify the proof, we consider $g$ and $h$ as functions $\mathbb{R}^3\to \mathbb{R}$, which makes no difference by the formulas of   $g$ and $h$. Also note that if a sequence approximates a proper line or a proper plane in $\mathbb{R}\times S^1\times S^1$, then its lift in $\mathbb{R}^3$ also approximates the lift of the proper line or  proper plane.

\noindent\textbf{Case 1.} Let the proper horizontal line be given by \eqref{eq:line}. Its image by $h$ is parametrized as 
$$h_1(t):=h(t,s^*,w^*)=\gamma^tB(2\pi w^*)\sin (2\pi s^*+\eta_2) - \xi(2\pi w^*).$$
It follows from the first inequality in \eqref{eq:con_line} that this image covers a neighborhood of zero. Since $cl(E)$ contains the proper horizontal line, by continuity of $h_1$, we obtain a subsequence for which  $\{h(t_n,s_n,w_n)\}$ is dense near zero.
Let us denote the subsequence by the same notations. By \eqref{eq:g} and \eqref{eq:h}, we have
\begin{equation}\label{eq:reduce}
\begin{aligned}
|g(t_n,s_n,w_n)|
&=\left|\dfrac{(h(t_n,s_n,w_n)+\xi(2\pi w_n))A(2\pi w_n)\sin(2\pi  s_n+\eta_1)}{B(2\pi w_n)\sin( 2\pi s_n+\eta_2)}\right|\\
&\to 
\left|\dfrac{\xi(2\pi w^*)A(2\pi w^*)\sin( 2\pi s^*+\eta_1)}{B(2\pi w^*)\sin( 2\pi s^*+\eta_2)}\right|.
\end{aligned}
\end{equation}
This gives the desired result on $|g(t_n,s_n,w_n)|$ due to the second inequality in \eqref{eq:con_line}.

\noindent\textbf{Case 2.} Let the proper tilted line be given by \eqref{eq:segment}. By \eqref{eq:short1}, one can find a small interval
 $\Delta_i$ arbitrarily close to $s=(-\eta_2/2\pi)+i$ with $i\in\mathbb{Z}$ such that 
\begin{equation}\label{eq:con_segment}
0<\xi(2\pi w)B(2\pi w)\sin(2\pi s+\eta_2)<\infty,\quad
\dfrac{\xi(2\pi w)A(2\pi w)\sin(2\pi s+\eta_1)}{\sin(2\pi s+\eta_2)}\not\in\{0,\pm 1,\infty\},
\end{equation}
for all $s\in \Delta_i$ with $w=a_2s$ or $w=w^*$.  Since $\xi(2\pi w)$ is  bounded, taking $n$ sufficiently large and $\Delta_i$ sufficiently close to $s=(-\eta_2/2\pi)+i$ ensures that the images of the functions
$$h_{21}(s):=h(a_1s+b,s,a_2 s)=\gamma^{a_1s+b}B(2\pi a_2s)\sin (2\pi s+\eta_2) - \xi(2\pi a_2s),\quad s\in\Delta_i$$
obtained with $w=a_2 s$ in \eqref{eq:segment}, and
$$h_{22}(s):=h(a_1s+b,s,w^*)=\gamma^{a_1s+b}B(2\pi w^*)\sin (2\pi s+\eta_2) - \xi(2\pi w^*),\quad s\in\Delta_i$$
obtained with $w=w^*$ in \eqref{eq:segment} both cover a small neighborhood of zero. Now, take a subsequence of $E$ approximating the line segment obtained by intersecting the tilted line with $\{(t,w,w)\mid s\in \Delta_i\}$. Obviously, for this subsequence $\{h(t_n,s_n,w_n)\}$ is dense near zero. Moreover, when $\Delta_i$ is taking sufficiently close, the value of $\sin(2\pi s+\eta_2)$ is arbitrarily close to zero, which together with \eqref{eq:reduce} implies that $\{|g(t_n,s_n,w_n)|\}$ lie in some closed interval in $(1,\infty)$.

\noindent\textbf{Case 3.} We now take a small interval $\Delta_i$ arbitrarily close to $s=(-\eta_2/2\pi) +i$ and a small interval $\Delta'_i$ arbitrarily close to $s=-\eta_1/2\pi +i$. Define two functions
$$h_{31}(t,s):=h(t,s,at+bs)=\gamma^{t}B(2\pi (at+bs))\sin (2\pi s+\eta_2) - \xi(2\pi(at+bs)),$$
in the first case of \eqref{eq:plane}, and
$$h_{32}(t,s):=h(t,s,w^*)=\gamma^{t}B(2\pi w^*)\sin (2\pi s+\eta_2) - \xi(2\pi w^*)$$
in the second case of \eqref{eq:plane}. Obviously, one can find $t^*$ such that, for all $s\in\Delta_i$ and either $w=at^*+bs$ or $w=w^*$, the inequalities in   \eqref{eq:con_segment} are satisfied and $h_{31}(\Delta_i)$ or $h_{32}(\Delta_i)$, depending on whether $w=at^*+bs$ or $w=w^*$, covers a small neighborhood of zero. Then, consider the subsequence approximating the intersection of the proper plane with $\{t=t^*, s\in\Delta_i\}$. Proceeding as in case 2 shows that $\{h(t_n,s_n,w_n)\}$ is dense near zero and $\{|g(t_n,s_n,w_n)|\}$ lie in some closed interval in $(1,\infty)$. Considering a similar number $t^*$ for $s\in\Delta_i'$ gives another required subsequence.
\end{proof}

\noindent{\it Proof of Proposition \ref{prop:blenderdf}.}
 Denote 
\begin{equation}\label{eq:sequence}
t_n=m_n-k_n\theta,\qquad s_n=k_n\dfrac{\omega_1}{2\pi} \bmod 1,\qquad w_n=m_n\dfrac{\omega_2}{2\pi} \bmod 1.
\end{equation}
Recall that in the double-focus case we have $\lambda>0,\gamma>0$, so $\lambda^k\gamma^m=\gamma^{m-k\theta}.$
Comparing $A_{k,m}$ and $B_{k,m}$ in \eqref{eq:df:AkmBkm} with functions $g$ and $h$ given by \eqref{eq:g} and \eqref{eq:h}, yields
$$A_{k_n,m_n}=g(t_n,s_n,w_n)
\quad\mbox{and}\quad
B_{k_n,m_n}=h(t_n,s_n,w_n).
$$
It follows that, for any sequence $\{(t_n,s_n,w_n)\}$ falling in one of the cases of Lemma~\ref{lem:reduction}, Proposition~\ref{prop:blender} is applicable to the system $\{T_n:=T_{k_n,m_n}\}$ (passing to the subsequence in the lemma) given by \eqref{eq:df:T_km_new}, with the correspondence $U\to X, W\to Y, V\to Z.$ More specifically, we have that
\begin{itemize}[nosep]
\item if $\{(t_n,s_n,w_n)\}$ approximates a proper horizontal line, then there exists a blender whose type depends on the coefficients in \eqref{eq:df:T_km_new};
\item if $\{(t_n,s_n,w_n)\}$ approximates a proper tilted  line, then there exists a cu-blender; and
\item if $\{(t_n,s_n,w_n)\}$ approximates a proper plane, then there exist simultaneously a cs-blender and a cu-blender.
\end{itemize}
Moreover, since the cone field in Proposition~\ref{prop:blender} does not depend on $\delta$ and $k$, one concludes (as in the proof of Proposition~\ref{prop:blendersf}) that the obtained blender has an activating pair $(\Pi',\mathcal{C}^{ss})$ if it is center-stable or $(\Pi',\mathcal{C}^{uu})$ if it is center-unstable, where $\Pi$ is defined in \eqref{eq:Pi'df} and $\mathcal{C}^{ss}$ and $\mathcal{C}^{uu}$ are given by Lemma~\ref{lem:conefields_df}.

We thus prove the Proposition \ref{prop:blenderdf} by showing that, for each of the cases in the statement of the proposition, one can find a sequence $\{(k_n,m_n)\}$ such that $\{(t_n,s_n,w_n)\}$ given by~\eqref{eq:sequence} approximates a proper line or  a proper plane. Recall that $\omega_{1,2}\in(0,\pi)$ and by $p/q$ we always mean that $p$ and $q$ are positive coprime integers.

\noindent\textbf{Case 1: $\theta$ is irrational, $\omega_1/2\pi,\omega_2/2\pi$ are rational.} 
Denote
$$\dfrac{\omega_1}{2\pi}=\dfrac{p_1}{q_1}\quad\mbox{and}\quad
\dfrac{\omega_2}{2\pi}=\dfrac{p_2}{q_2}.
$$
Since condition \ref{word:A2} is satisfied, there exist integers $r_1\in\{0,\dots,q_1-1\}$ and $r_2\in\{0,\dots,q_2-1\}$ such that inequalities in \eqref{eq:con_line} are satisfied with $s^*=r_1{\omega_1}/{2\pi}$ and $w^*=r_2\omega_2/2\pi$. Then, 
$$\left\{s=r_1\dfrac{\omega_1}{2\pi}\bmod 1,w=r_2\dfrac{\omega_2}{2\pi}\bmod 1\right\}$$
is a proper horizontal line of the form \eqref{eq:line}. The irrationality of $\theta$ implies that we can take a sequence $\{(k'_n,m'_n)\}$ of pairs of positive integers such that the set
$$\{q_2m'_n+r_2-(q_1k'_n+r_1)\theta\}$$
is dense in $\mathbb{R}$. Setting $k_n=q_1k'_n+r_1$ and $ m_n=q_2m'_n+r_2$ in \eqref{eq:sequence} gives a sequence whose closure contains the above proper horizontal line. So, we get a cs-blender if $|\alpha_{r_1,r_2}|<1$, or a cu-blender if $|\alpha_{r_1,r_2}|>1$.

\noindent\textbf{Case 2.1: $\omega_1/2\pi$ is irrational, $\omega_2/2\pi$ is rational, and $\theta,\omega_1/2\pi,1$ are rationally independent.} Let $p_2,q_2,r_2$ be as in the preceding case. By the rational independence among $\theta,\omega_1/2\pi,1$, one can find a sequence $\{(k'_n,m'_n,i'_n)\in\mathbb{N}^3\}$ such that the set
$$\left\{\left(q_2m'_n+r_2-k'_n\theta,k'_n\dfrac{\omega_1}{2\pi} - i'_n\right)\right\}$$
dense in $\mathbb R^2$. Setting $k_n=k'_n$ and $m_n=q_2m'_n+r_2$ in \eqref{eq:sequence} gives a sequence whose closure contains the proper plane $\{w=r_2\omega_2/2\pi\}.$ So, we get simultaneously a cs-blender and a cu-blender.\\

Before proceed to other cases, we make the following 
\begin{claim*}
To prove that there exists $\{(t_n,s_n,w_n)\}$ given by \eqref{eq:sequence} approximating a proper line or a proper plane, it suffices to show that each point of this proper line or proper plane is accumulated by points of the set
\begin{equation}\label{eq:testset}
\left\{\left(m_n-k_n\theta,k_n\dfrac{\omega_1}{2\pi}\bmod 1,m_n\dfrac{\omega_2}{2\pi}\bmod 1\right)\right\}
\end{equation}
for some sequence $\{(k_n,m_n)\}$.
\end{claim*}
\noindent Indeed, the assumption particularly implies that, for each rational point in the proper line or proper plane, we have such a sequence $\{(k_n,m_n)\}$, and hence the union of these countably many sequences gives a new sequence which defines the desired  $\{(t_n,s_n,w_n)\}$ by \eqref{eq:sequence}.\\

\noindent\textbf{Case  2.2: $\omega_1/2\pi$ is irrational, $\omega_2/2\pi$ is rational, and $\theta,\omega_1/2\pi,1$ are rationally dependent.} Denote
\begin{equation}\label{eq:df:2.2r}
\theta =\dfrac{p_1}{q_1} \dfrac{\omega_1}{2\pi} + \dfrac{p_2}{q_2}\quad\mbox{and}\quad
\dfrac{\omega_2}{2\pi}=\dfrac{p_3}{q_3}.
\end{equation}
With $r_2$  as in case 1, the set $\{t=-sp_1/q_1+r_2,w=r_2p_3/q_3\bmod 1\}$ is a proper tilted line. Let us show that each point of this line can be approximated by the  set \eqref{eq:testset}.

By the irrationality of $\omega_1/2\pi$, the set
\begin{equation*}
\left\{q_2q_3k\dfrac{\omega_1}{2\pi}-q_1q_3i\right\}_{k,i\in\mathbb N^2}
\end{equation*}
is dense in $\mathbb{R}$. So, for any $s\in [0,1]$, one can find a sequence $\{(k_n',i_n')\}$ such that 
$$q_2q_3k'_n\dfrac{\omega_1}{2\pi}-q_1q_3i'_n\to s.$$
Denote $k_n=q_2q_3k'_n $ and $ m_n=r_2+p_2q_3k'_n+p_1q_3i'_n.$ The relation \eqref{eq:df:2.2r} then implies
\begin{align*}
m_n-k_n\theta &=(r_2+p_2q_3k'_n+p_1q_3i'_n)-\dfrac{p_1}{q_1} k_n \dfrac{\omega_1}{2\pi} - p_2q_3k_n'\\
&=r_2 -\dfrac{p_1}{q_1}\left(k_n \dfrac{\omega_1}{2\pi}-q_1q_3i'_n\right)\\
&\to r_2-s\dfrac{p_1}{q_1},\\[10pt]
m_n\dfrac{\omega_2}{2\pi}\mod 1&=(r_2+p_2q_3k'_n+p_1q_3 i'_n)\dfrac{p_3}{q_3} \mod 1\\
&=\dfrac{r_2p_3}{q_3} \mod 1,
\end{align*}
which gives the desired approximation. By the claim, we have a cu-blender.

\noindent\textbf{Case 3.1: $\omega_1/2\pi,\omega_2/2\pi$ are irrational, and $\theta,\omega_1/2\pi,\theta\omega_2/2\pi,1$ are rationally independent.}
The rational independence implies that, for any point $(t,s,w)$, there exists $\{(k_n,m_n,i_n,\ell_n)\}$ such that
\begin{equation}\label{eq:conver2}
m_n-k_n\theta\to t,\quad
 k_n\dfrac{\omega_1}{2\pi}-i_n\to s,\quad
k_n\theta\dfrac{\omega_2}{2\pi}-\ell_n\to w- t\dfrac{\omega_2}{2\pi}.
\end{equation}
Combining the first and third relations, yields
$$(m_n-t)\dfrac{\omega_2}{2\pi}-\ell_n \to w-t\dfrac{\omega_2}{2\pi},$$
which implies 
\begin{equation}\label{eq:conver1}
m_n-k_n\theta\to t,\quad
 k_n\dfrac{\omega_1}{2\pi}-i_n\to s,\quad
m_n\dfrac{\omega_2}{2\pi}-\ell_n\to w,
\end{equation}
Therefore, any point of any proper plane can be approximated by a sequence of the form \eqref{eq:testset}. By the claim, we get simultaneously a cs-blender and a cu-blender.

\noindent\textbf{Case 3.2: $\omega_1/2\pi,\omega_2/2\pi$ are irrational, $\theta,\omega_1/2\pi,\theta\omega_2/2\pi,1$ are rationally dependent, but $\theta,\omega_1/2\pi,1$ are rationally independent.} Denote
\begin{equation}\label{eq:df:3.2r}
\theta\dfrac{\omega_2}{2\pi}=
\dfrac{p_1}{q_1}\theta + \dfrac{p_2}{q_2}\dfrac{\omega_1}{2\pi}+\dfrac{p_3}{q_3}.
\end{equation}
Since $\theta,\omega_1/2\pi,1$ are rationally independent, the set
$$\left\{\left(q_1m-q_3k\theta,q_3k\dfrac{\omega_1}{2\pi}-q_2 i\right)\right\}_{k,m,i\in\mathbb{N}}$$
is dense in $\mathbb{R}^2$. Hence, for any point $(t,s)$ there exists a sequence $(k'_n,m'_n,i'_n)$ such that 
\begin{equation*}
q_1m'_n-q_3k'_n\theta\to t \quad\mbox{and}\quad
 q_3k'_n\dfrac{\omega_1}{2\pi}-q_2i'_n \mod 1\to s.
\end{equation*}
Denote  $k_n=q_3k'_n $ and $ m_n=q_1m'_n.$ The relation \eqref{eq:df:3.2r} shows that
\begin{align*}
 m_n\dfrac{\omega_2}{2\pi}\mod 1
&= k_n\left(\dfrac{p_1}{q_1}\theta + \dfrac{p_2}{q_2}\dfrac{\omega_1}{2\pi}+\dfrac{p_3}{q_3}\right)-\dfrac{\omega_2}{2\pi}t\mod 1\\
&\to \left(\dfrac{p_1}{q_1}-\dfrac{\omega_2}{2\pi}\right)t+\dfrac{p_2}{q_2}s\mod 1\\
&=:a t+bs\mod 1.
\end{align*}
Thus, any point of the proper plane $\{w=at+bs\bmod 1\}$ can be approximated by a sequence of the form \eqref{eq:testset}. By the claim, we get simultaneously a cs-blender and a cu-blender.

\noindent\textbf{Case 3.3.1: $\omega_1/2\pi,\omega_2/2\pi$ are irrational, $\theta,\omega_1/2\pi,1$ are rationally dependent, and $\omega_1/2\pi,\theta\omega_2/2\pi,1$ are rationally dependent.} Denote
\begin{equation}\label{eq:df:3.3.1r}
\theta=\dfrac{p_1}{q_1}\dfrac{\omega_1}{2\pi}+\dfrac{p_2}{q_2}
\quad\mbox{and}\quad
\theta\dfrac{\omega_2}{2\pi}= \dfrac{p_3}{q_3}\dfrac{\omega_1}{2\pi}+\dfrac{p_4}{q_4}.
\end{equation}
The irrationality of $\omega_1/2\pi$ implies that the set
$$\left\{q_2q_4k\dfrac{\omega_1}{2\pi}-q_1q_3i\right\}_{k,i\in\mathbb N}$$
is dense in $\mathbb{R}$. For any $s\in [0,1]$ and any sequence $(k'_n,i'_n)$ such that 
\begin{equation*}
q_2q_4k'_n\dfrac{\omega_1}{2\pi}-q_1q_3i'_n \to s,
\end{equation*}
denote $k_n=q_2q_4k'_n $ and $ m_n=p_2q_4k'_n+p_1q_3i'_n.$
The relations in \eqref{eq:df:3.3.1r} show that
\begin{align*}
m_n-k_n\theta
&=p_2q_4k'_n+p_1q_3i'_n-q_2q_4k'_n\dfrac{p_1}{q_1}\left(\dfrac{\omega_1}{2\pi}+\dfrac{p_2}{q_2}\right)\\
&=-\dfrac{p_1}{q_1}\left(q_2q_4k'_n\dfrac{\omega_1}{2\pi}-q_1q_3i'_n\right)\\
&\to -\dfrac{p_1}{q_1} s,\\[10pt]
k_n\theta\dfrac{\omega_2}{2\pi}\mod 1
&=\dfrac{p_3}{q_3}\left(q_2q_4k'_n\dfrac{\omega_1}{2\pi}-q_1q_3i'_n\right)+p_3q_1i'_n+p_4q_2k'_n\mod 1\\
&\to \dfrac{p_3}{q_3}s \mod1.
\end{align*}
Observing the relation between  \eqref{eq:conver2} and \eqref{eq:conver1}, one finds
$$m_n\dfrac{\omega_2}{2\pi}\to \left(\dfrac{p_3}{q_3}-\dfrac{\omega_2}{2\pi}\dfrac{p_1}{q_1}\right)s=:as.$$
Thus, any point of the proper tilted line $\{t=-p_1s/q_1,w=as\}$ can be approximated by a sequence of the form \eqref{eq:testset}.
The claim then gives  a cu-blender.

\noindent\textbf{Case 3.3.2: $\omega_1/2\pi,\omega_2/2\pi$ are irrational, $\theta,\omega_1/2\pi,1$ are rationally dependent, and $\omega_1/2\pi,\theta\omega_2/2\pi,1$ are rationally independent.} We still denote $\theta$ as in \eqref{eq:df:3.3.1r}. For any $s\in [0,1]$  and any sequence $(k'_n,i'_n)$ such that
$$q_2k_n'\dfrac{\omega_1}{2\pi}-q_1i'_n\to s,$$
denote $k_n=q_2k_n' $ and $ m_n=p_2k'_n+p_1i'_n.$
The first relation in \eqref{eq:df:3.3.1r} gives
$$m_n-k_n\theta 
=p_2k'_n+p_1i'_n-
q_2k_n'\left(
\dfrac{p_1}{q_1}\dfrac{\omega_1}{2\pi}+\dfrac{p_2}{q_2}\right)
\to - \dfrac{p_1}{q_1}s.$$
It then follows from the rational independence among $\omega_1/2\pi,\theta\omega_2/2\pi,1$ and the correspondence relation between \eqref{eq:conver2} and \eqref{eq:conver1} that any point of any proper tilted line in the plane $\{t=-p_1s/q_1\}$ can be approximated by a sequence of the form  \eqref{eq:testset}. The claim then gives a cu-blender.
\qed

\section{Robust heterodimensional dynamics in generic one-parameter unfolding families: proof of Theorem \ref{thm:unfolding1}}\label{sec:unfolding}
Finding robust heterodimensional dynamics is essentially equivalent to connecting the invariant manifolds of the blenders and periodic orbits  $O_{i}$ ($i=1,2$). As shown in Sections \ref{sec:homosf} and \ref{sec:homodf}, the structure of the blenders obtained in this paper  is quite clear in the small neighborhood $\Pi$ of $M^+_1$ (for saddle-focus cycles) and of $M^-_2$ (for double-focus cycles), with explicit formulas for their local invariant manifolds. The homoclinic relations or activations between blenders and the  points $O_{i}$ depend on how the local invariant manifolds of $O_i$ behave along the dynamics of the cycle. By construction, the intersections of $W^s(O_1)$ with the unstable manifolds of the blenders and the intersections of $W^u(O_2)$ with the stable manifolds of the blenders are automatic (see \eqref{eq:autoint1} and \eqref{eq:autoint2}). The other intersections depend on the positions of the iterates 
$$F_{21}\circ F^m_2\circ F_{12}(W^u_{loc}(O_1))
\quad\mbox{and}\quad
 F^{-k}_1\circ F^{-1}_{12}(W^s_{loc}(O_2))$$
for saddle-focus cycles, and
$$F^{m}_2\circ F_{12}(W^u_{loc}(O_1))
\quad\mbox{and}\quad
 F^{-1}_{21}\circ F^{-k}_1\circ F^{-1}_{12}(W^s_{loc}(O_2))$$
for double-focus cycles. Basically, one gets the intersection of $W^u(O_1)$ with the stable manifolds of the blenders if the first iterate in each case crosses $\Pi'$, and gets  the intersection of $W^s(O_2)$ with the unstable manifolds of the blenders if the second iterate in each case crosses $\Pi'$. According to Lemma~\ref{lem:LTlemmasf}, the second iterate in the saddle-focus case always hits $\Pi'$ if $\omega/2\pi$ is irrational. Similarly, in the double-focus case, Lemmas~\ref{lem:LTlemma14} and~\ref{lem:LTlemma15}  show that the two iterates cross $\Pi'$ if $\omega_2/2\pi$ and, respectively, $\omega_1/2\pi$ is irrational. Thus, the proof of Theorem~\ref{thm:rational2} is essentially about verifying the intersections of these iterates with $\Pi'$ when $\omega/2\pi$ and $\omega_i/2\pi$ are rational and the cycle is unfolded.

Recall that $\mu$ is the splitting parameter of a non-degenerate heterodimensional cycle, see Remark~\ref{rem:mu0}.  A one-parameter unfolding family $\{f_\eps\}$ with $f_0=f$ having a non-degenerate cycle satisfies $d\mu(0)/d\eps\neq 0$ (see~\eqref{eq:rank}). Since it suffices to prove the results for one-parameter families, we simply use $\mu$ as a new parameter, replacing  $\eps$. 

\subsection{Unfolding saddle-focus heterodimensional cycles}\label{sec:unfoldsf}
Let us first deal with saddle-focus cycles, where the center-stable multipliers are nonreal and the center-unstable one is real, denoted by $\lambda e^{i\omega}$ and $\gamma$, respectively.  
 For brevity, we denote below the continuations of hyperbolic objects (e.g. $O_{1,2}$ and the blenders) by the same letters and drop the term `continuation'.

Let the coefficient $A,B,u^-,b,\eta_{1,2}$ be as in \eqref{akmbkm}, let the cube $\Pi'$ be given by \eqref{eq:Pi'sf}, and let the cone fields $\mathcal{C}^{ss}$ and $\mathcal{C}^{uu}$ be given by Lemma~\ref{lem:conefields}. 

\begin{lem}[{\citep[Lemma 6.4]{LT21}}\footnote{As mentioned, the coordinates in \citep{LT21} are $(X_1,X_2,Y,Z)$, and in this paper we dropped the subscript of $X_1$ and grouped $X_2,Z$ as a new $Z$ coordinate. Moreover, the coefficient $b_{11}$ there is denoted by $b$ in this paper.}]\label{lem:unfoldO1}
Fix any sufficiently small $\delta$. For every sufficiently large $m$, if
\begin{equation}\label{eq:unfold1:1}
|\mu - \gamma^{-m}u^- |<\dfrac{1}{2}|\gamma^{-m}b^{-1}|q\delta,
\end{equation}
then $F_{21}\circ F_2 \circ F_{12}(W^u_{loc}(O_1))$ is a disc of the form $(X,Z)=s(Y)$ for some smooth function $s$, which crosses $\Pi'$ properly with respect to $\mathcal{C}^{uu}$.
\end{lem}

\begin{lem}\label{lem:unfoldO2}
Suppose $\omega/2\pi=p/q$. Fix a sufficiently small $\delta$ and any $r\in \{0,\dots,q-1\}$. For every $k=r+iq$ with  sufficiently large $i$, if
\begin{equation}\label{eq:unfold2:1}
|\mu + \lambda^{k}b^{-1}B\sin(2\pi r p/q+\eta_2)|
<\dfrac{1}{2}|\lambda^kA\sin(2\pi r p/q+\eta_1)|q\delta,
\end{equation}
then $F^{-k}_1\circ F^{-1}_{12}(W^s_{loc}(O_2))$ is a disc of the form $(X,Y)=s(Z)$ for some smooth function $s$, which crosses $\Pi'$  properly with respect to $\mathcal{C}^{ss}$. 
\end{lem}

The proof of Lemma~\ref{lem:unfoldO2}  is given in Appendix~\ref{sec:lem:unfoldO2} as a minor modification of the proof of \citep[Lemma 6.5]{LT21}. 

\begin{proof}[Proof of Theorem \ref{thm:unfolding1} for the saddle-focus case.] The results in Theorem \ref{thm:unfolding1} for the saddle-focus case are those concerning Theorems \ref{thm:sf} and \ref{thm:rational2}.
By \eqref{eq:rank}, we only need to consider the family $\{f_\mu\}$ with $f_0=f$. 
Since blenders persist under $C^1$-small perturbations, the homoclinic relations and activations obtained in Theorems \ref{thm:sf} and \ref{thm:rational2} are preserved for all small $\mu$.

\noindent\textbf{Case (1) of Theorem \ref{thm:sf}.} It suffices to show that $O_1$ activates the cu-blender $\Lambda^{cu}$ and is homoclinically related to some non-trivial hyperbolic basic set. Then, by the lambda-lemma this set also activates $\Lambda^{cu}$. 
Since $W^s_{loc}(O_1)=\{Y=0\}$ by \eqref{eq:Omf}, Lemma \ref{lem:unfoldO1} implies that there exists a sequence $\mu_j\to 0$ such that $F_{21}\circ F_2 \circ F_{12}(W^u_{loc}(O_1))\cap W^s_{loc}(O_1)\neq \emptyset$, which gives a hyperbolic basic set homoclinically related to $O_1$. To see that $O_1$ actives $\Lambda^{cu}$, we note first that $W^s(O_1)\cap W^u(\Lambda^{cu})\neq\emptyset$ by \eqref{eq:Omf} and the formula for $W^u_{loc}(M_2)$ in Section \ref{sec:homosf}. Then, since $\Lambda^{cu}$ is obtained from Proposition \ref{prop:blendersf} and hence has an activating pair $(\Pi',\mathcal{C}^{uu})$, 
the proper crossing part of Lemma \ref{lem:unfoldO1} shows that $O_1$ indeed activates $\Lambda^{cu}$.

\noindent\textbf{Case (2) of Theorem \ref{thm:sf}.} One only needs to create a homoclinic relation between the cs-blender $\Lambda^{cs}$ and $O_1$ by changing $\mu$. The intersection $W^s(O_1)\cap W^u(\Lambda^{cs})\neq\emptyset$ is automatic at $\mu=0$ by \eqref{eq:Omf} and the formula for $W^u_{loc}(M_1)$ in Section \ref{sec:homosf}. With the formula for $W^s_{loc}(M_1)$ in the same section, one obtains by Lemma \ref{lem:unfoldO1} a sequence $\mu_j\to 0$ for which $W^u(O_1)\cap W^s(\Lambda^{cs})\neq \emptyset.$

\noindent\textbf{Case (3) of Theorem \ref{thm:sf} and Theorem~\ref{thm:rational2}.}
Before proceeding to the rational $\omega$ cases, let us make an observation. By similar arguments to the above cases, we obtain that, when Lemma \ref{lem:unfoldO2} holds, $O_2$ is homoclinically related to the cu-blender in Proposition \ref{prop:blendersf}, if exists; it activates the cs-blender, if exists; and, it is always homoclinically related to a non-trivial hyperbolic basic set. We have thus proved that\\ \\
\textit{
for every sufficiently small $\mu$ satisfying both conditions \eqref{eq:unfold1:1} and \eqref{eq:unfold2:1},  each of $O_1$ and $O_2$ is homoclinically related to a non-trivial hyperbolic basic set; and, moreover, whenever there is a cs-blender  from Proposition \ref{prop:blendersf}, it is homoclinically related to $O_1$ and activated by $O_2$, and whenever there is a cu-blender from Proposition \ref{prop:blendersf}, it is homoclinically related to $O_2$ and activated by $O_1$.}\\ \\
Then, the results concerning Theorem \ref{thm:rational2} and case (3) of Theorem \ref{thm:sf} follow immediately from
\begin{claim}\label{lem:mu}
There exists a sequence $\mu_j\to 0$ with each $\mu_j$ satisfying both conditions \eqref{eq:unfold1:1} and \eqref{eq:unfold2:1}.
\end{claim}
This completes the proof of Theorem \ref{thm:unfolding1} for saddle-focus  cycles. Let us prove the above claim. Note that~\eqref{eq:unfold1:1} and~\eqref{eq:unfold2:1} give two intervals $I^u_m$ and, respectively, $I^s_k$ of $\mu$ values such that the inequality is satisfied if $\mu$ lies in corresponding interval. It then suffices to find a sequence $\{k_n,m_n\}$ with $k_n,m_n\to\infty$ such that $I^u_{m_n}\cap I^s_{k_n}\neq \emptyset$. The centers of these two intervals are given by
$$
c^u_m = \gamma^{-m}u^- \quad\mbox{and}\quad c^s_k=-\lambda^kb^{-1}B\sin(k\omega+\eta_2),
$$
and the radius  of $I^u_m$ is $|\gamma^{-m}b^{-1}|q\delta/2$. Now set $k_n=k'_nq+r$ and $m_n=m'_n$ as in \eqref{eq:sf:sequence}. It follows from~\eqref{eq:lambdagamma} that
$$
\dfrac{|c^u_m-c^u_k|}{|\gamma^{-m}b^{-1}|q\delta/2}
=2 \cdot \dfrac{|\lambda^{k_n}\gamma^{m_n}b^{-1}B\sin(k_n\omega+\eta_2)-u^-|}{b^{-1}\delta'} \to 0 \quad\mbox{as}\quad n\to\infty,
$$
which shows that $I^u_{m_n}\cap I^s_{k_n}\neq \emptyset$ for all sufficiently large $n$.
\end{proof}

\subsection{Unfolding double-focus heterodimensional cycles}
Here we have the following  counterpart  to  Lemmas \ref{lem:unfoldO1} and \ref{lem:unfoldO2}, and Claim \ref{lem:mu}:

\begin{lem}
Let the cube $\Pi'$ be given by \eqref{eq:Pi'df}, and the cone fields $\mathcal{C}^{uu}$ and $\mathcal{C}^{ss}$ be given by Lemma~\ref{lem:conefields_df}.
\begin{itemize}[nosep]
\item If $\omega_2/\pi$ is rational, then there exists a sequence of intervals $I^u_m$ such that, for all sufficiently large $m$, if $\mu\in I^u_m$ then $F^{m}_2\circ F_{12}(W^u_{loc}(O_1))\cap \Pi'$ is a disc of the form $(U,V)=s(W),$
crossing $\Pi'$ properly with respect to $\mathcal{C}^{uu}$.
\item  If $\omega_1/\pi$ is rational, then there exists a sequence of intervals $I^s_k$ such that, for all sufficiently large $k$, if $\mu\in I^s_k$ then $F^{-1}_{21}\circ F^{-{k}}_1\circ F^{-1}_{12}(W^s_{loc}(O_2))\cap \Pi'$ is a disc of the form $(U,W)=s_2(V),$
crossing $\Pi'$ properly with respect to $\mathcal{C}^{ss}$.
\item If both $\omega_1/\pi$ and $\omega_2/\pi$ are rational, then there exists $\{k_n,m_n\}$ such that $I^u_{m_n}\cap I^s_{k_n}\neq \emptyset$ for all large $n$.
\end{itemize}
\end{lem}
The first two items  are perturbative versions of Lemma \ref{lem:LTlemma14} and  Lemma \ref{lem:LTlemma15}, respectively. They can be obtained by modifying the proofs in \citep[Section 6.4.5]{LT21} in the same way as we did to obtain  Lemma \ref{lem:unfoldO2}, the perturbative version of Lemma \ref{lem:LTlemmasf}. The third item can also be proved in the same way as Claim \ref{lem:mu}. We hence omit their proofs.  Theorem \ref{thm:unfolding1} for double-focus cycles then follows from the routine checks as in Section \ref{sec:unfoldsf}, with using the above results and the lemmas in Section \ref{sec:homodf}. \qed


\section{Purely rational moduli: proof of Theorem \ref{thm:rational}}\label{sec:rational}

Since the set of non-degenerate cycles is open and dense in $\mathcal{H}_f$ (see Appendix~\ref{sec:nd}), we can always assume that the cycle in Theorem \ref{thm:rational} is non-degenerate. The first item of Theorem \ref{thm:rational} concerns the saddle case, i.e., when both the center-stable multiplier $\lambda$ and center-unstable multiplier $\gamma$ are real, which follows immediately from {\citep[Theorem 6]{LT21}}.

%

We now consider the saddle-focus and double-focus cases. We need to investigate the set of orbits which lies entirely in a small neighborhood $U$ of a heterodimensional cycle. Recall that a cycle is the union $\Gamma=L_1\cup L_2\cup \Gamma^0\cup\Gamma^1$, where $L_{1,2}$ are the orbits of the periodic points $O_{1,2}$ and $\Gamma^0$ and $\Gamma^1$ are the two heteroclinic orbits in $W^u(L_1)\cap W^s(L_2)$ and $W^u(L_2)\cap W^s(L_1)$, respectively. An orbit which lies entirely in $U$, except for those in the invariant manifolds of $L_{1,2}$, must intersect the neighborhood $\Pi$ of $M^+$ infinitely many times for both forward and backward iterations. In particular, an intersection of $\Pi$ with any orbit gives a pair of two points $M$ and $\bar M$ in this orbit satisfying $\bar M=T_{k,m}(M)$. In what follows we show that (Lemma \ref{lem:allrat}), for any sequence of pairs $\{(k_i,m_i)\}$ corresponding to intersection points of $\Pi$ with any set of orbits, the assumption of rational moduli imposes a strong relation on the pairs, which thus leads to simple hyperbolic dynamics.

\subsection{Simple hyperbolic dynamics in the saddle-focus case}
Let us first deal with the saddle-focus cycles, where the center-stable and center-unstable multipliers are $\lambda e^{\pm i\omega}$ and $\gamma$, and we assume that $\theta$ and $\omega/2\pi$ are rational.  

Consider a generic unfolding family $\{f_\mu\}$. To study the dynamics, we use the following perturbative version of formula \eqref{eq:complex:T_k,m_cross_0mu}, which is easily derived from \citep{LT21} in Appendix~\ref{sec:formula}:
any two points $(X,Y,Z),(\bar X,\bar Y,\bar Z)\in\Pi$
(the domain defined in \eqref{eq:complex:domain}) satisfy $(\bar X,\bar Y,\bar Z)=T_{k,m}(X,Y,Z)$ for some pair $(k,m)$ if and only if 
\begin{equation}\label{eq:revision1}
\begin{aligned}
\bar X &= b\gamma^{m}\mu+  A_{k,m} X + B_{k,m} +\phi_{1}(X,\bar Y,Z),\\
 Y &=\phi_2(X,\bar Y,Z),\quad \bar Z = \phi_3(X, \bar Y,Z),
\end{aligned}
\end{equation}
where the coefficients and the functions are the same as in \eqref{eq:complex:T_k,m_cross_0mu},  with smooth dependence on $\mu$, and  $\phi_1$ satisfies the finer estimate
\begin{equation}\label{eq:phi1}
\phi_{1}=\gamma^m(o(\lambda^k)+o(\gamma^{-m})).
\end{equation}

Let $\omega/2\pi=p/q$. There exist two finite sets $E_{\eta_i}\subset \mathbb R$ ($i=1,2$) such that if $\eta_{i}\not\in E_{\eta_i}$, then for any $r \in \{0,\dots, q-1\}$ one has
\begin{equation}\label{eq:etacon}
\sin(2\pi r p/q+ \eta_1)\neq 0
\quad\mbox{and}\quad
\sin(2\pi r p/q+ \eta_2)\neq 0.
\end{equation}
For these  $\eta_i\notin E_{\eta_i}$, we use \eqref{akmbkm} to rewrite the above formula as
\begin{equation}\label{eq:sfmu}
\begin{aligned}
\bar X &= b\gamma^{m}\mu+  a_rb\lambda^k\gamma^m   X + a_rbx^+_r\lambda^k\gamma^m- b u^- +\phi_{1}(X,\bar Y,Z),\\
 Y &=\phi_3(X,\bar Y,Z),\quad \bar Z = \phi_4(X, \bar Y,Z),
\end{aligned}
\end{equation}
where
\begin{equation}\label{akmbkm2}
a_r=b^{-1}A \sin(2\pi r p/q+ \eta_1),\quad
x^+_r=\dfrac{B\sin(2\pi r p/q+ \eta_2)}{A \sin(2\pi r p/q+ \eta_1)},\quad r\in\{0,\dots,q-1\}.
\end{equation}

\begin{lem}\label{lem:allrat}
Let $\theta=p'/q'$ and $\omega/2\pi=p/q$. Assume that for all $r_1,r_2\in\{0,\dots,q-1\}$ we have
\begin{equation}\label{eq:allratcon}
\left|\dfrac{u^-}{a_{r_1} x^+_{r_1}}\right|\not\in \bigcup_{i=0}^{q-1}cl\left\{ 
|\gamma|^{\frac{s}{q'}}\cdot\dfrac{|{a_{r_2} x^+_{r_2}}(a_{r_1}{x^+_{r_1}})^{-1}|-\lambda^\ell}{1-\gamma^{-n}}
 \right\}_{s\in\mathbb{Z},\ell\in\mathbb{N},n\in \mathbb{N}}=:E_{r_1,r_2}.
\end{equation}
Consider any (finite or infinite) sequence of different pairs $\{(k_i,m_i)\}$ such that there exist some sequence of points $\{M_i\in \Pi\}$ and some small $\mu$ for which $T_{k_i,m_i}(M_i)\in\Pi$ with $\Pi$  defined in \eqref{eq:complex:domain}. Then, for all sufficiently small $\delta$, either all $m_i$ are equal to each other, $m_i=m$, and $\lambda^{k_i} = O(\delta) \gamma^{-m} $ for all $i$, or all $k_i$ are equal to each other, $k_i=k$, and $\gamma^{-m_i} = O(\delta) \lambda^{k}$ for all $i$.
\end{lem}
\begin{proof}
Let us first consider the special case where the sequence $\{(k_i,m_i)\}$ has two pairs. By assumption the system 
\begin{equation}\label{eq:twopt}
\bar M_1=T_{k_1,m_1}(M_1),\qquad \bar M_2=T_{k_2,m_2}(M_2)
\end{equation}
have at least one solution $(k_1,m_1,M_1,\bar M_1,k_2,m_2,M_2,\bar M_2)$ with $k_i,m_i\in \mathbb{N}$ and $M_i,\bar M_i\in \Pi$. This system imposes two equations involving the $X$-equation in  \eqref{eq:sfmu}, which are
\begin{equation*}
\begin{aligned}
b\gamma^{m}\mu+  a_{r_1}b\lambda^{k_1}\gamma^{m_1}   X_1 + a_{r_1}bx^+_{r_1}\lambda^{k_1}\gamma^m- b u^- +\gamma^m(o(\lambda^k)+o(\gamma^{-m}))&=\bar X_1,\\
b\gamma^{m}\mu+  a_{r_2}b\lambda^{k_2}\gamma^{m_2}   X_2 + a_{r_2}bx^+_{r_2}\lambda^{k_2}\gamma^{m_2}- b u^- +\gamma^m(o(\lambda^k)+o(\gamma^{-m}))&=\bar X_2,
\end{aligned}
\end{equation*}
for some $r_1,r_2\in\{0,\dots,q-1\}$. Since $X_i,\bar X_i \in [-\delta,\delta] $ by assumption, the new system has solutions if and only if there exist $X_i=K_i\delta$ and $\bar X_i=C_i\delta$ for some $K_i,C_i\in[-1,1]$ (where $\delta$ is the size of $\Pi$ defined in \eqref{eq:complex:domain}) such that 
\begin{equation*}
\begin{aligned}
\mu+ a_{r_1}  \lambda^{k_1}(x^+_{r_1} + K_1\delta ) - \gamma^{-m_1} u^- &= \gamma^{-m_1}b^{-1}C_1\delta,\\
\mu+ a_{r_2} \lambda^{k_2}(x^+_{r_2} + K_2\delta ) - \gamma^{-m_2} u^- &= \gamma^{-m_2}b^{-1}C_2\delta.
\end{aligned}
\end{equation*}
Subtracting the second equation for the first one yields
\begin{equation}\label{eq:twopt2}
\lambda^{k_1}-
\lambda^{k_2}\dfrac{a_{r_2}(x^+_{r_2} + K_2\delta )}{a_{r_1} (x^+_{r_1} + K_1\delta )}
=
\gamma^{-m_1}\dfrac{C_1\delta+bu^-}{a_{r_1}b(x^+_{r_1} + K_1\delta )}
-
\gamma^{-m_2}\dfrac{C_2\delta+bu^-}{a_{r_1}b(x^+_{r_1} + K_1\delta )}.
\end{equation}

From this point, the proof is the same as that  of \citep[Lemma 5.7]{LT21}. Indeed, comparing \eqref{eq:twopt2} with (5.36) in \citep{LT21}, we see that the only difference is that $a_{r_1}=a_{r_2}=a$ and $x^+_{r_1}=x^+_{r_2}=x^+$ in \citep{LT21}. Modifying the proof there  in the obvious way, one finds that equation  \eqref{eq:twopt2} has no solutions satisfying $(k_1,m_1)\neq (k_2,m_2)$  unless $$
\left|\dfrac{u^-}{a_{r_1} x^+_{r_1}}\right|\not\in cl\left\{ 
|\gamma|^{\frac{s}{q'}}\cdot\dfrac{|{a_{r_2} x^+_{r_2}}(a_{r_1}{x^+_{r_1}})^{-1}|-\lambda^\ell}{1-\gamma^{-n}}
 \right\}_{s\in\mathbb{Z},\ell\in\mathbb{N},n\in \mathbb{N}}.
$$
Proceeding as in   \citep{LT21} with considering all possible pairs $(r_1,r_2)$, we obtain the desired conclusion on $k_i$ and $m_i$ for an arbitrary two-element sequence $\{(k_1,m_1),(k_2,m_2)\}$. 

To finish the proof of the lemma, let us consider a sequence $\{(k_i,m_i)\}$ with more than two elements. Apply the above result to the first two elements, we have $m_1=m_2=m$ or $k_1=k_2=k$. We assume the former case and the proof in the latter case is the same. Then, we have
$$
\lambda^{k_1}=O(\delta)\gamma^{-m_1}.
$$
We claim that $m_i=m$ for all $i$. If this were not true, then there would be some $i$ such that $m_i\neq m=m_1$. Applying the two-element result to $\{(k_1,m_1),(k_i,m_i)\}$, yields $k_1=k_i$ and
$\gamma^{-m_1}=O(\delta)\lambda^{k_1}$, contradicting the above inequality. Thus, all $m_i$ are the same, and the relation $\lambda^{k_i}=O(\delta)\gamma^{-m}$ is obtained by applying the two-element result  to $\{(k_i,m),(k_j,m)\}$ for every $(i,j)$ pair.
\end{proof}

Recall that $E_{\eta_{i}}$ ($i=1,2$) are the set of $\eta_i$ values for which \eqref{eq:etacon} is satisfied, and $E_{r_1,r_2}$ are the sets in \eqref{eq:allratcon}.
\begin{thms}\label{thm:appendix}
Let $f$ have a non-degenerate saddle-focus heterodimensional cycle with both $\theta=p'/q'$ and $\omega/2\pi=p/q$ being rational. If
\begin{equation*}
\eta_{i}\notin E_{\eta_{i}},\quad
\left|\dfrac{u^-}{a_{r_1} x^+_{r_1}}\right|\notin E_{r_1,r_2},\quad
|a_{r}b|\not\in\{|\gamma|^{\frac{s}{q'}}\}_{s\in\mathbb{Z}} =: E'_{r},
\end{equation*}
with $i\in\{1,2\}$ and  $r_1,r_2,r\in\{0,\dots,q-1\}$, then any generic one-parameter unfolding family $\{f_\eps\}$ with $f_0=f$ has simple hyperbolic dynamics as in Definition \ref{defi:simple}.
\end{thms}
Since each of the sets $E_{\eta_i}$, $E_{r_1,r_2}$ and $E'_r$ is nowhere dense in $\mathbb{R}$, this theorem covers  the second and third items of Theorem \ref{thm:rational}.  The proof is essentially a word-by-word repetition of the proof of \citep[Theorem 6]{LT21}, with formula \eqref{eq:sfmu} in place of \citep[equation (3.20)]{LT21} and   Lemma \ref{lem:allrat} in place of \citep[Lemma 5.7]{LT21}. 
For completeness, we give the proof in Appendix~\ref{sec:rationalapp}.

\subsection{Simple hyperbolic dynamics in the double-focus case} 
The arguments for double-focus cycles are completely parallel. In this case, the central multipliers are  $\lambda e^{\pm i\omega_1}$ and $\gamma e^{\pm i\omega_2}$, and  
$$\theta=\dfrac{p'}{q'},\quad
\dfrac{\omega_1}{2\pi}=\dfrac{p_1}{q_1},\quad
\dfrac{\omega_2}{2\pi}=\dfrac{p_2}{q_2}
$$
are rational. First note that there is a set $E_{\eta_3}$ such that if $\eta_3\in E_{\eta_3}$, then the condition \eqref{eq:quant} is satisfied. For these $\eta_3$ values, one can obtain a perturbative version of formula \eqref{eq:df:T_km_new} by performing the calculation in the proof of \citep[Lemma 6.6]{LT21} with adding the parameter $\mu$ to the right hand side of the first equation in (6.22) there. We obtain that, for point $(U,V,W)\in \Pi$ (defined in \eqref{eq:df:domain}), one has 
\begin{equation}\label{eq:twopt3}
\begin{array}{l}
 \bar U 
= \ell_{r_2}\lambda^k\mu + A_{k,m}U+B_{k,m} + \phi_1(U,V,\bar W),\\
\bar{V}=\phi_2(U,V,\bar W),\qquad
W=\phi_3(U,V,\bar W),
\end{array}
\end{equation}
where $A_{k,m},B_{k,m}$ and the functions $\phi$ are as in \eqref{eq:df:T_km_new} with smooth dependence on $\mu$, and 
$$\ell_{r_2}=\dfrac{2}{\sqrt{1+\ell^2}\cos(2\pi r_2 p_2/q_2+\eta_3)},\quad r_2\in\{0,\dots,q_2-1\}.$$
with $\ell$
being certain coefficient\footnote{The coefficient`$a_{14}$' in equation (6.22) of \citep{LT21}}  in $\D F_{12}$. Denoting 
\begin{align*}
&a_{r_1,r_2}= A(2\pi r_2 p_2/q_2)\sin(2\pi r_1 p_1/q_1+\eta_1),\\
&b_{r_1,r_2}=B(2\pi r_2 p_2/q_2)\sin(2\pi r_1 p_1/q_1+\eta_2),\qquad
c_{r_2}=\xi(2\pi r_2 p_2/q_2),
\end{align*}
with $r_i\in \{0,\dots,q_i - 1\}$ ($i=1,2$), one rewrites \eqref{eq:twopt3} as
\begin{equation*}\label{eq:twopt4}
\begin{array}{l}
 \bar U 
= \ell_{r_2}\lambda^k\mu +  a_{r_1,r_2}\lambda^k\gamma^mU+  b_{r_1,r_2}\lambda^k\gamma^m +c_{r_2} + \phi_1(U,V,\bar W),\\
\bar{V}=\phi_2(U,V,\bar W),\qquad
W=\phi_3(U,V,\bar W).
\end{array}
\end{equation*}
Arguing as in the saddle-focus case, one  easily obtains a counterpart to Lemma \ref{lem:allrat}, and hence obtains the following result covering the last item of  Theorem~\ref{thm:rational}:

\begin{thms}
Let $f$ have a non-degenerate double-focus heterodimensional cycle with $\theta=p'/q',\omega_1/2\pi=p_1/q_1,\omega_2/2\pi=p_2/q_2$ being rational. There exist nowhere dense subsets of $\mathbb{R}$: $E_{\eta_i},E_{r_1,r_2,r'_1,r'_2},E'_{r_1,r_2}$ with $i\in\{1,2,3\}$, $r_1,r'_1\in \{0,\dots,q_1-1\}$ and $r_2,r'_2\in\{0,\dots,q_2-1\}$ such that if
\begin{equation*}
\eta_{i}\notin E_{\eta_{i}},\quad
\left|\dfrac{1}{b_{r_1,r_2}}\right|\notin E_{r_1,r_2,r'_1,r'_2},\quad
|a_{r_1,r_2}|\notin E'_{r_1,r_2},
\end{equation*}
then any generic one-parameter unfolding family $\{f_\eps\}$ with $f_0=f$ has simple hyperbolic dynamics as in Definition \ref{defi:simple}.
\end{thms}

\section*{Acknowledgment}
\addcontentsline{toc}{section}{Acknowledgments}
I am grateful to Professor Dmitry Turaev for his invaluable guidance during my career and for bringing me the interesting topics on heterodimensional cycles.  I  thank Xiaolong Li for useful discussions that improved the presentation of the main theorems. I also thank the anonymous referee for  helpful comments and suggestions. This work was supported by the Science Fund Program for Excellent Young Scientists (Overseas), the New Cornerstone Science Foundation, and  the Leverhulme Trust grant RPG 2021-072.

\appendix
\section{Non-degeneracy conditions}\label{sec:nd}

 We use the notations in Section~\ref{sec:setting}.  Let us first collect some preliminary facts on the local dynamics of hyperbolic periodic points (see e.g. \citep{SSTC1}). In $U_{01}$, the local map $F_1$ has a $(d-d_1)$-dimensional local stable invariant manifold $W^s_{loc}(O_1)$ and a $d_1$-dimensional local unstable manifold $W^u_{loc}(O_1)$. On $W^s_{loc}(O_1)$ there exists a strong-stable $C^r$-smooth foliation $\mathcal{F}^{ss}$, which particularly includes the $(d-d_1-2)$-dimensional strong-stable invariant manifold $W^{ss}_{loc}(O_1)$. There also exists a $(d_1+2)$-dimensional extended-unstable invariant manifold $W^{uE}_{loc}(O_1)$ tangent at $O_1$ to the eigenspace  corresponding to the two conjugate nonreal center-stable multipliers and  the unstable multipliers $\gamma_{1,1},\dots ,\gamma_{1,d_1}$ (recall that the multipliers are introduced in \eqref{eq:mult} and in the saddle-focus case we always assume that $L_1$ has nonreal center-stable multipliers). The manifold $W^{uE}_{loc}(O_1)$  is not uniquely defined; however, any two such manifolds contain $W^u_{loc}(O_1)$ and are tangent to each other at points of $W^u_{loc}(O_1)$.  

Similarly, there exists a strong-unstable $C^r$-smooth foliation $\mathcal{F}^{uu}$ on the $d_2$-dimensional local unstable invariant manifold $W^u_{loc}(O_2)$. Denote by $d^{cu}$ the number of center-unstable multipliers of $O_2$, so $d^{cu}=1$ if the cycle is of saddle-focus type, or $d^{cu}=2$ if it is of double-focus type. The manifold $W^u_{loc}(O_2)$ includes the $(d_2-d^{cu})$-dimensional strong-unstable invariant manifold $W^{uu}_{loc}(O_2)$. Here we also have a $(d-d_2+d^{cu})$-dimensional extended-stable invariant manifold $W^{sE}_{loc}(O_2)$, corresponding to the center-unstable multiplier(s)  and the stable multipliers $\lambda_{2,1},\dots ,\lambda_{2,d-d_2}$. Any two such manifolds contain $W^s_{loc}(O_2)$ and are tangent to each other at the points of $W^s_{loc}(O_2)$.

We introduce the following non-degeneracy conditions:

\noindent\setword{\textbf{C1}}{word:C1}: $F^{-1}_{12}(W_{loc}^{sE}(O_2))\pitchfork W^{u}_{loc}(O_1)\neq\emptyset$ at the point $M^-_1$ and
$F_{12}(W_{loc}^{uE}(O_1))\pitchfork W^{s}_{loc}(O_2)\neq\emptyset$ at $M^+_2$;\\[5pt]
\noindent\setword{\textbf{C2}}{word:C2}: the leaf of $\mathcal{F}^{uu}$ at the point $M^-_2$ is not tangent to
$F_{21}^{-1}(W^s_{loc}(O_1))$ and the leaf of $\mathcal{F}^{ss}$ at the point $M^+_1$ is not tangent to
$F_{21}(W^u_{loc}(O_2))$, and \\[5pt]
\noindent\setword{\textbf{C3}}{word:C3}: $\Gamma^1\cap (W^{ss}(O_1)\cup W^{uu}(O_2))=\emptyset$, i.e., the points $M^+_1$ and $M^-_2$ lie outside of $W^{ss}(O_1)$ and $W^{uu}(O_2)$, respectively. 

Figure~\ref{fig:gc} provides an illustration of the fulfilment of these three conditions.
 The first two conditions are about the non-degeneracy of the non-transverse and transverse intersections. Condition \ref{word:C2} is $C^r$-open, and condition \ref{word:C1} is $C^r$-open in the $\mathcal{H}_f$,  the class of diffeomorphisms with the heterodimensional cycle defined above Theorem \ref{thm:rational}. One can always make $f$ fulfill conditions \ref{word:C1} and \ref{word:C2} by an arbitrarily small $C^r$-perturbation; more specifically, in the smooth case one adds a local perturbation to $f$ and in the analytic case one follows the scheme used in \citep{BT86,GTS07}. Condition \ref{word:C3} is also $C^r$-open, and, if it is not satisfied, then, without perturbing $f$, we can always find another transverse heteroclinic orbit which satisfies \ref{word:C3} and is arbitrarily close to $\Gamma^1$. To see this, note that  condition \ref{word:C2} shows that the curve $\ell_1=W^s_{loc}(O_1)\cap F_{21}(W^u_{loc}(O_2))$ is not tangent to $W^{ss}_{loc}(O_1)\cup F_{21}(W^{uu}_{loc}(O_2))$. Hence, there is always another point in $\ell_1$ curve lying outside of that $W^{ss}(O_1)$ and $W^{uu}(O_2)$. We then take this orbit of this point as the new transverse heteroclinic orbit in our cycle.

There is one more non-degeneracy condition. To state it precisely, let us introduce $C^r$-coordinates $(x,y,z)\in \mathbb{R}^{2}\times\mathbb{R}^{d_1}\times\mathbb{R}^{d-d_1-2}$ in $U_{01}$ such that the local  invariant manifolds as well as the leaves of the foliation $\mathcal{F}^{ss}$ are straightened (see Figure \ref{fig:gc}). So we have
$$W^s_{loc}(O_1)=\{y=0\},\quad W^u_{loc}(O_1)=\{x=0, z=0\},\quad W^{ss}_{loc}(O_1)=\{x=0, y=0\}.$$
These coordinates can be chosen such that the extended-unstable manifold $W^{uE}_{loc}(O_1)$ is tangent to $\{z=0\}$ when $x=0$ and $z=0$.

Similarly, we introduce $C^r$-coordinates 
$(u,v,w)\in\mathbb{R}^{d^{cu}}\times\mathbb{R}^{d-d_2}\times\mathbb{R}^{d_2-d^{cu}}$ in $U_{02}$ such that the local  invariant manifolds as well as the leaves of the foliation $\mathcal{F}^{uu}$ are straightened, so we have
$$W^s_{loc}(O_2)=\{u=0, w=0\},\quad W^u_{loc}(O_2)=\{v=0\},\quad W^{uu}_{loc}(O_2)=\{u=0,v=0\}.$$
These coordinates can be chosen such that the extended-stable manifold $W^{sE}_{loc}(O_2)$ is tangent to $\{w=0\}$ when $u=0$ and $v=0$.

With  \cite[Lemmas 5 and 6]{GST08}, the coordinates can be  chosen to further satisfy  that the restrictions
$$F_1|_{W^s_{loc}(O_1)}:(x,z)\mapsto (\bar x,\bar z)\quad\mbox{and}\quad F_2|_{W^u_{loc}(O_2)}:(u,w)\mapsto (\bar u,\bar w)$$ are linear in the center-stable and, respectively, center-unstable coordinates:
\begin{equation}\label{eq:intro:2}\!\!\!\!
\begin{array}{ll}
\begin{array}{c}\mbox{Saddle-focus}\\\mbox{case:}\end{array}&\!\bar{x} = \lambda\begin{pmatrix}
\cos \omega &\sin\omega \\
-\sin \omega & \cos\omega
\end{pmatrix}\!x \quad\mathrm{and}\quad \bar{u}=\gamma u;\\[15pt]
\begin{array}{c}\mbox{Double-focus}\\\mbox{case:}\end{array}&\!\!\bar{x} = \lambda\!\begin{pmatrix}
\!\cos \omega_1 &\!\!\sin\omega_1 \\
\!-\sin \omega_1 &\!\!\cos\omega_1
\end{pmatrix}\!x \;\;\mathrm{and}\;\; \bar{u}=\gamma\!\begin{pmatrix}
\!\cos \omega_2 &\!\!\sin\omega_2 \\
\!-\sin \omega_2 &\!\!\cos\omega_2
\end{pmatrix}\! u.\!\!\!\!\!\!\!\!\!\!\!
\end{array}
\end{equation}
Recall that for  saddle-focus cycles we omit the case where the center-stable multiplier is real but the center-unstable multipliers are nonreal, since it can be reduced to the above case by considering $f^{-1}$. In these coordinates, we have $M^+_1=(x^+,0,z^+)$ and $M^-_2=(u^-,0,w^-)$, and hence condition \ref{word:C3} reads $x^+\neq 0$ and $u^-\neq 0$.  

Recall that the transition map $F_{21}$ takes a small neighborhood of $M^-_2$ to a small neighborhood of $M^+_1$.   By counting dimension, the two intersections $\ell_1=W^s_{loc}(O_1)\cap F_{21}(W^u_{loc}(O_2))$ and $\ell_2=F^{-1}_{21}\ell_1$ are two one-dimensional curves. Denote by $(dx,dz)_{M^+_1}$ the tangent vectors  to $\ell_1$ at the point $M^+_1$, and by $(du,dw)_{M^-_2}$ the tangent vectors to $\ell_2$ at the point $M^-_2$.

\noindent\setword{\textbf{C4}}{word:C4}: when the center-stable multipliers are a pair of complex conjugate numbers (so $x\in\mathbb{R}^2$), the vector $dx_{M^+_1}$ is not parallel to  the vector $x^+$. When the center-unstable multipliers are a pair of complex conjugate numbers (so $u\in\mathbb{R}^2$), the vector $du_{M^-_2}$ is not parallel to the vector $u^-$.


This condition is independent of the choice of the coordinates which straighten the local invariant manifolds and foliations, and  linearize the central dynamics as in \eqref{eq:intro:2}. To see this, one notices that any transformation between  two sets of such coordinates leaves the $x$- and $u$-plane invariant and restricts to a linear map on them.
Condition \ref{word:C4} is also independent of the choice of the points $M^+_1$ and $M^-_2$ from $\Gamma^1$, since replacing them by their iterates contribute the same factor to $dx_{M^+_1}$ and  $x^+$, and, respectively, to  $du_{M^-_2}$ and  $u^-$.

A saddle-focus or double-focus heterodimensional cycle is called {\em non-degenerate} if it satisfies conditions C1--C4, and any such cycle can be made non-degenerate by an arbitrarily $C^r$-small perturbation.

\begin{rem}\label{rem:gcsaddle}
For heterodimensional cycles of saddle type, i.e., both center-stable and center-unstable multipliers are real, the first three non-degeneracy conditions are the same. The fourth condition is different. In this case, the central dynamics on $W^s_{loc}(O_1)$ and $W^u_{loc}(O_2)$ are given by $\bar{x} = \lambda x,\,\bar{u}=\gamma u$ with $x\in \mathbb{R}$ and $u\in \mathbb{R}$. By condition \ref{word:C2}, the restriction $F_{21}|_{\ell_2}:\ell_2\to\ell_1,(u,v,w)\mapsto(x,v,z)$ has the following non-zero derivative:
\begin{equation*}\label{eq:intro:b}
b^*= \left.\frac{\partial x}{\partial u}\right|_{M^-_2}.
\end{equation*}
The fourth condition is $\left|{b^*u^-}/{x^+}\right|\neq 1$ (which is condition C4.1 in \citep[Section 2.3]{LT21}). 
\end{rem}


\section{Simple results analogous  to those in \citep{LT21}}
\subsection{Sketch of the proof of Proposition~\ref{prop:blender}}\label{sec:blendercon}
Let us first prove the proposition when condition 2.1 is satisfied. The map $T_n^\times: (X,\bar Y,Z)\mapsto (\bar X,Y,\bar Z)$ given by \eqref{eq:blender:1} is a contraction. Take a finite set $\mathcal{N}\subset \mathbb N$ such that $T_n^\times$ takes $Q$ into itself  for all $n\in \mathcal{N}$. Then, by a fixed point lemma for a direct product of metric spaces (see e.g. Theorem 6.2 in \citep{Sh67}), there is a one-to-one correspondence between the set $\Sigma:=\{0,\dots,\card\mathcal{N}-1\}^\mathbb{Z}$ of two-sided sequences of $\card\mathcal{N}$ symbols and the set $\hat \Lambda$ of points whose orbits under the iterated function system $\{T_n\}_{n\in \mathcal{N}}$ never leave $Q$. This shows that  $\hat \Lambda$ is a compact, transitive, locally maximal invariant set of $T$. Evidently from  formula \eqref{eq:blender:1}, $\hat \Lambda$ is also hyperbolic.

Let $(\Delta X,\Delta Y, \Delta Z)$ denote vectors in tangent spaces. Given any constant $K>0$, it is easy to check that, by taking $\delta$ sufficiently small and $n$ sufficiently large, the cone field defined by
\begin{equation}\label{eq:ssconeblender}
\mathcal{C}^{ss}=\{(\Delta X,\Delta Y, \Delta Z): |\Delta X|+\|\Delta Y\|<K\|Z\|\},
\end{equation}
is backward invariant in the sense that for any $(X,Y,Z)\in\Pi$ if $T_n(X,Y,Z)=(\bar X,\bar Y,\bar Z)\in\Pi$ for some $n\in\mathcal{N}$, then the the cone at $(\bar X,\bar Y,\bar Z)$ is mapped inside the one at $( X, Y, Z)$ by $\D T^{-1}_n$. 

By assumption 1, there exists $c\in(0,1)$ such that $B_n$ is dense in $[-c\delta,c\delta]$. Let $Q'$ be defined as in \eqref{eq:Q}.  Consider the open set $\mathcal{S}$ of discs $S$ of the form $(X,Y)=s(Z)$ for some smooth function $s$, which crosses $Q'$ properly with respect to $\mathcal{C}^{ss}$.
We will prove that the set $\Lambda=\hat\Lambda\cap Q'$ is a cs-blender, which amounts to showing that $S\cap W^u(\Lambda)\neq\emptyset$ for every $S\in \mathcal{S}$. By construction, the local unstable manifold of $\Lambda$ is the set of points in $Q$ whose backward orbits under $\{T_n\}_{n\in \mathcal{N}}$ stay in $Q'$. Thus, it suffices to prove the so-called {\em covering property} that for any $S\in \mathcal{S}$ there exists $n\in\mathcal{N}$ such that the intersection $T^{-1}_n(S)\cap Q'$ contains a piece in $\mathcal{S}$. Indeed, this gives a sequence of positive integers $\{n_i\in\mathcal{N}\}$ and a sequence of nested compact sets $\{\hat S_i\subset S\}$ such that each preimage  $T_{n_i}(\hat S_i)$ is a disc in $\mathcal{S}$, which particularly implies that $\bigcap \hat S_i$ contains a point whose backward orbit under $\{T_n\}_{n\in \mathcal{N}}$ lies in $Q$.

Let us prove the covering property. Due to the backward invariance of $\mathcal{C}^{ss}$, one only needs to check that $T^{-1}_n(S)\cap Q$ has a connected component crossing $Q$. Also observe that by the strong contraction in $Z$ and strong expansion in $Y$, it suffices to show that, for any $S\in\mathcal{S}$, there exists $n\in\mathcal{N}$ such that  the $X$-coordinates of points of $T^{-1}_n(S)$ lie in $I:=[-c\delta,c\delta]$. This is equivalent to prove that
\begin{itemize}[nosep]
\item for $n\in\mathcal{N}$, the union of the images of the interior of $I$ under the affine maps
$g_n:\bar X = A_{n} X + B_{n}$
covers a neighborhood of $I$, and that
\item the overlaps between $g_n(I)$ are large enough so that there exists some $n^*\in \mathcal{N}$ such that the $X$-coordinates of points of $S$ lie completely in $g_{n^*}(I)$.
\end{itemize}
To find a set $\mathcal{N}$ satisfying the above conditions, we note that the center of $g_n(I)$ is $B_n$ and the length of it is greater than some constant $L_1>0$. Since $\{B_n\}$ is dense in $I$, it follows that there exist finitely many $n$ such that the the union of the $L_1/2$-neighborhoods of $B_n$ covers a neighborhood of $I$. The overlap requirement can be satisfied by taking $K$ in \eqref{eq:ssconeblender} sufficiently small so that the total variation of $X$-coordinates of $S$ is sufficiently small. The covering property is  proven.

Finally, one sees that the blender found above can be arbitrarily close to the origin by taking the integers in $\mathcal{N}$ sufficiently large. This completes the proof when condition 2.1 holds.

When condition 3.1 holds, define $\hat T_n:= T^{-1}_n|_{T_n(\sigma_n)\cap Q}$. The maps $\hat T_n$ assume the same form as \eqref{eq:blender:1} but with new $A_n$ whose absolute values are uniformly bounded away from 0 and 1. Then, apply the above arguments to $\{\hat T_n\}$.

\subsection{Proof of Lemma~\ref{lem:unfoldO2}}\label{sec:lem:unfoldO2}
In the  coordinates \eqref{eq:intro:2}, we have $M^-_1=(0,y^-,0)$ and $M^+_2=(0,v^+,0)$ for some constants $y^-\neq 0$ and $v^+\neq 0$. The parameter $\mu$ enters the Taylor expansion of the transition map $F_{12}$ as
\begin{equation*}
\begin{aligned}
u &= \mu + a_{11} x + a_{12}(y-y^-) + a_{13} z+\dots \\
v-v^+ &=  a_{21} x + a_{22}(y-y^-) + a_{23} z+\dots \\
w &=  a_{21} x + a_{32}(y-y^-) + a_{33} z+\dots, 
\end{aligned}
\end{equation*}
where all coefficients $y^-,v^+,a_{ij}$ depend on $\mu$.  
\begin{rem}\label{rem:mu}
This is formula (6.3) of \citep{LT21}, where the parameter is denoted by $\hat{\mu}$. When $r\geq 3$ (which is our case), it is the same as our $\mu$, see \citep[Section 6.1.2]{LT21}.
\end{rem}

We need a formula for $F^{-k}_1\circ F^{-1}_{12}(W^s_{loc}(O_2))$. Such formula at $\mu=0$ is  given by  (6.30) of \citep{LT21}, and a parametric version with non-zero $\mu$ can be obtained by setting $\hat\mu$ (which is the $\mu$ here by Remark~\ref{rem:mu}) non-zero in (6.8) of \citep{LT21} and then applying the coordinate transformation (6.11) of \citep{LT21}. In the end, one finds that $F^{-k}_1\circ F^{-1}_{12}(W^s_{loc}(O_2))$ is given by
\begin{equation*}\label{eq:unfold2:2}
\begin{aligned}
X &=\dfrac{-b\lambda^{-k}\mu- B\sin(2\pi r p/q+\eta_2)}{A\sin(2\pi r p/q+\eta_1)} + O(\delta^{\frac{3}{2}}) + o(1)_{k\to\infty}, \\ 
Y&=o(\hat\gamma^{-k}),
\end{aligned}
\end{equation*}
where $\hat{\gamma}>|\gamma|$ is some constant and the right-hand sides are functions of $Z$ defined on $[-\delta,\delta]^{d-d_1-1}$. It immediately follows that, for all sufficiently small $\delta$ and sufficiently large $k$, if $\mu$ satisfies \eqref{eq:unfold2:1}, then 
$$|X|<q\delta,\quad \|Y\|<\delta,\quad \dfrac{\partial (X,Y)}{\partial Z}=O(\delta^{\frac{3}{2}}) + o(1)_{k\to\infty},$$
where the first two inequalities mean that $F^{-k}_1\circ F^{-1}_{12}(W^s_{loc}(O_2))$ crosses $\Pi'$ and the last one implies that the crossing is proper with respect to $\mathcal{C}^{ss}$.

\subsection{Remarks on formula \eqref{eq:revision1}}\label{sec:formula}
This formula can be obtained  by first applying the coordinate transformation (6.11) in \citep{LT21} to equation (6.10) there with keeping the $\hat \mu$ term, then making the notation change $b_{11}\to b,X_1\to X, (X_2,Z)\to Z$, and finally using the Implicit Function Theorem in an obvious way to express the nonlinearities as functions of $(X,\bar Y,Z)$. The estimate \eqref{eq:phi1}  is given by equation (6.12) in \citep{LT21}.

One notices that here we do not have condition \eqref{eq:crosscon}. Indeed, this condition is equivalent to condition (6.17) in \citep{LT21}, which is  to ensure that the image $T_{k,m}(M)$ of a point $M\in\Pi$ lies also in $\Pi$, and this fact then allows one to use the Implicit Function Theorem to deduce formula (6.18) (i.e., formula \eqref{eq:complex:T_k,m_cross_0mu} in the present paper) form (6.12). Since we now directly assume   $\bar M=(\bar X,\bar Y,\bar Z)\in \Pi$,  condition \eqref{eq:crosscon} is no longer needed.

\subsection{Proof of Theorem~\ref{thm:appendix}}\label{sec:rationalapp}
First note that Lemma \ref{lem:allrat} holds when $k_i$ or $m_i$ is infinite, with taking $\lambda^\infty=0$ and $\gamma^{-\infty}=0$. The case $k_i=\infty$ corresponds to $M_i\in W^s_{loc}(O_1)$ and $\bar M_i \in F_{21}\circ F_2^{m_i} \circ F_{12} (W^u_{loc}(O_1))$.
The case $m_i=\infty$ corresponds to $M_i \in F_1^{-k_i}\circ F_{12}^{-1} (W^s_{loc}(O_2))$ and $\bar M_i \in F_{21}(W^u_{loc}(O_2))$.

The case $k_1=m_1=\infty$ corresponds to the existence of the  heteroclinic intersection of $F_{12}(W^u_{loc}(O_1))$ with $W^s_{loc}(O_2)$, where one can take $M_1\in W^s_{loc}(O_1)$ and $\bar M_1\in T_{21}(W^u_{loc}(O_2))$. This happens at $\mu=0$. Suppose there is another pair $(k_2,m_2)$ such that   $\bar M_2=T_{k_2,m_2}(M_2)$. Then, one can easily check that applying Lemma \ref{lem:allrat} to the sequence $\{(k_1,m_1),(k_2,m_2)\}$ forces $k_2=m_2=\infty$, which implies the dynamics described in as in Definition \ref{defi:simple} for $\mu=0$.

For $\mu\neq 0$, consider the set $\mathcal{L}$ of orbits which lie entirely in $U$ and different from $L_{1,2}$ and from the heteroclinic orbits in the intersection of $F_{21}(W^u_{loc}(O_2))$ with $W^s_{loc}(O_1)$. Any orbit $L\in \mathcal{L}$ gives a sequence $\{M_i\}$ of intersection points with $\Pi$, satisfying $M_{i+1}=T_{k_i,m_i}(M_i)$ for some $(k_i,m_i)$. Note that the case $k_i=m_i=\infty$ is excluded since $F_{12}(W^u_{loc}(O_1))\cap W^s_{loc}(O_2)=\emptyset$ at $\mu\neq 0$. Depending on $\{M_i\}$, the sequence $\{(k_i,m_i)\}$ can be infinite or finite. We  modify the sequence $\{(k_i,m_i)\}$ as follows such that it is always infinite.
\begin{itemize}[topsep=0pt]
\item If $\{M_i\}$ is infinite, then the sequence $\{(k_i,m_i)\}$ is infinite with finite $k_i,m_i$, and we keep it;
\item If $\{M_i\}$ is finite from the left, then the leftmost point $M_{i_l}$ lies either in $F_{21}(W^u_{loc}(O_2))$ or in $F_{21}\circ F_2^{m} \circ F_{12} (W^u_{loc}(O_1))$ for some $m$. But the former is impossible since $L\in\mathcal{L}$. In this case, we extend $\{(k_i,m_i)\}$ to an infinite sequence by adding $(k_{i_l-j},m_{i_l-j})=(\infty,m)$ for all $j\geq 1$;
\item If $\{M_i\}$ is finite from the left, then the rightmost point $M_{i_r}$ lies either in $W^s_{loc}(O_1)$ or in $F_1^{-k}\circ F_{12}^{-1} (W^s_{loc}(O_2))$ for some $k$. But the former is impossible since $L\in\mathcal{L}$. In this case, we extend $\{(k_i,m_i)\}$  by adding $(k_{i_l+j},m_{i_l+j})=(k,\infty)$ for all $j\geq 0$.
\end{itemize}
Lemma \ref{lem:allrat} then implies three possibilities:
\begin{enumerate}[topsep=0pt]
\item there exist positive integers $k$ and $m$ such that $k_i=k$ and $m_i=m$ for all $i$ and for all $L\in\mathcal{L}$. In this case, one achieves the desired result by using the fact that the map \eqref{eq:sfmu} is hyperbolic when $|a_rb\lambda^k\gamma^m|$ is bounded away from 0 and 1, which  is equivalent to requiring $|a_rb|\not\in\{|\gamma|^{\frac{s}{q'}}\}$, and hence the maps $T_{k,m}$ are hyperbolic with contraction in $(X,Z)$ and expansion in $Y$ when $0<|a_rb\lambda^k\gamma^m|<1$, and with contraction in $Z$ and expansion in $(X,Y)$ when $1<|a_rb\lambda^k\gamma^m|$.
\item there exists positive integer $m$ such that $m_i=m$ for all $i$ and for all $L\in\mathcal{L}$, while some $k_i$ are different. The desired result follows from the fact that $|a_rb\lambda^k\gamma^m|<1$ by the relation $\lambda^{k_i}=O(\delta)\gamma^{-m}$ given by Lemma \ref{lem:allrat}, and hence the maps $T_{k_i,m}$ are hyperbolic with contraction in $(X,Z)$ and expansion in $Y$;
\item There exists a positive integer $k$ such that $k_i=k$ for all $i$ and for all $L\in\mathcal{L}$, while some $m_i$ are different. The desired result follows from the fact that $|a_rb\lambda^k\gamma^m|>1$ by the relation $\gamma^{-m_i}=O(\delta)\lambda^{k}$, and hence the maps $T_{k_i,m}$ are hyperbolic with contraction in $Z$ and expansion in $(X,Y)$.
\end{enumerate}

\end{document}